\documentclass[trsc,nonblindrev]{informs3} 

\OneAndAHalfSpacedXI 

\usepackage{natbib}
 \bibpunct[, ]{(}{)}{,}{a}{}{,}%
 \def\bibfont{\small}%

\usepackage{lineno,hyperref}
\usepackage{multirow}
\usepackage{etex}
\usepackage{amsfonts}
\usepackage{mathrsfs}
\usepackage{comment}
\usepackage{verbatim} 
\usepackage{tablefootnote}
\usepackage{float}
\usepackage{algorithm}
\usepackage{algpseudocode}
\newcommand{\algstep}[1]{\item[]\medskip\hrule\kern 2pt\hbox to \textwidth{\hspace{\labelsep}\textbf{#1}\hfill}\hrule}
\usepackage{todonotes}
\usepackage{amsmath}
\usepackage{longtable}
\usepackage{threeparttable}
\usepackage{caption}
\usepackage{subfigure}

\usepackage{xspace}
\newcommand\BRFrag{battery-restricted fragment\xspace}
\newcommand\BRFrags{battery-restricted fragments\xspace}

\def\TheoremsNumberedThrough{%
    \theoremstyle{TH}%
    \newtheorem{theorem}{Theorem}
    
    \newtheorem{proposition}{Proposition}

    \theoremstyle{EX}
    
    \newtheorem{exmp}{Example}
    
    \newtheorem{definition}{Definition}

}

\usepackage{microtype}
\usepackage{amsmath}
\usepackage{amssymb}
\usepackage{graphicx}
\usepackage{stmaryrd}
\usepackage{color}
\usepackage{footnote}
\makesavenoteenv{table}
\newcommand{\beq}{\begin{equation}}
\newcommand{\eeq}{\end{equation}}



\makeatletter 
\def\@makecaption#1#2{%
  \vskip\abovecaptionskip
  \sbox\@tempboxa{#1 #2}%
    {\bfseries #1} #2\par
  \vskip\belowcaptionskip}
\makeatother
\usepackage{booktabs}
\usepackage[justification=centering]{caption}

\usepackage{hyperref}
\hypersetup{
    colorlinks=true,
    linkcolor=blue,
    filecolor=gray,      
    urlcolor=blue,
    citecolor=blue,
}

\TheoremsNumberedThrough     

\EquationsNumberedThrough    

\begin{document}


\RUNAUTHOR{Su, Puchinger, and Dupin}

\RUNTITLE{A DA Algorithm for the E-ADARP}

\TITLE{A Deterministic Annealing Local Search for the Electric Autonomous Dial-A-Ride Problem}

\ARTICLEAUTHORS{%
\AUTHOR{Yue Su}
\AFF{Université Paris-Saclay, CentraleSupélec, Laboratoire Génie Industriel, 91190, Gif-sur-Yvette, France, \EMAIL{yue.su@centralesupelec.fr}}

\AUTHOR{Jakob Puchinger}
\AFF{EM Normandie Business School, Métis Lab, 92110, Clichy, France, \EMAIL{jpuchinger@em-normandie.fr}, \URL{}}
\AFF{Université Paris-Saclay, CentraleSupélec, Laboratoire Génie Industriel, 91190, Gif-sur-Yvette, France, \EMAIL{jakob.puchinger@centralesupelec.fr}}

\AUTHOR{Nicolas Dupin}
\AFF{Univ Angers, LERIA, SFR MATHSTIC, F-49000 Angers, France, \EMAIL{nicolas.dupin@univ-angers.fr}}
} 

\ABSTRACT{
This paper investigates the Electric Autonomous Dial-A-Ride Problem (E-ADARP), which consists in designing a set of minimum-cost routes that accommodates all customer requests for a fleet of Electric Autonomous Vehicles (EAVs). Problem-specific features of the E-ADARP include: (i) the employment of EAVs and a partial recharging policy; (ii) the weighted-sum objective function that minimizes the total travel time and the total excess user ride time. We propose a Deterministic Annealing (DA) algorithm, which is the first metaheuristic approach to solve the E-ADARP. Partial recharging (i) is handled by an exact route evaluation scheme of linear time complexity. To tackle (ii), we propose a new method that allows effective computations of minimum excess user ride time by introducing a fragment-based representation of paths. These two methods compose an exact and efficient optimization of excess user ride time for a generated E-ADARP route. To validate the performance of the DA algorithm, we compare our algorithm results to the best-reported Branch-and-Cut (B\&C) algorithm results on existing instances. Our algorithm provides 25 new best solutions and 45 equal solutions on 84 existing instances. To test the algorithm performance on larger-sized instances, we establish new instances with up to 8 vehicles and 96 requests, and we provide 19 new solutions for these instances. Our final investigation extends the state-of-the-art model and explores the effect of allowing multiple visits to recharging stations. This relaxation can efficiently improve the solution's feasibility and quality.

}%


\KEYWORDS{dial-a-ride problem; electric autonomous vehicles; deterministic annealing, metaheuristic}

\maketitle


%

\section{Introduction}
With the astounding growth of automobile ownership, a series of transport-related problems has appeared worldwide. These problems, such as greenhouse gas emissions and urban traffic congestion, have severely impacted the economy and the environment \citep{schrank20122012}. 
One possible approach to address these concerns is to provide ride-sharing services \citep{jin2018ridesourcing}, which require customers to specify their origins and destinations. 
The underlying optimization problem is usually modeled as a Dial-A-Ride Problem (DARP), which consists in designing minimum-cost routes for a fleet of vehicles to serve a set of customer requests \citep{cordeau2007dial}. Each customer request contains an origin, a destination, and a time window on either the origin or the destination. The DARP was first introduced in \cite{wilson1971scheduling} and has received considerable attention from the literature  \citep{parragh2008survey,molenbruch2017typology,ho2018survey}. The standard version of the DARP aims to minimize the total routing cost while respecting operational constraints such as time windows, capacity, and duration constraints. However, as customers can share rides with others, user inconvenience must be considered while minimizing the total routing cost. In the typical DARP model, a maximum user ride time constraint is introduced for each customer request. Due to the integration of maximum user ride time and time window constraints, scheduling vehicles to begin their services as early as possible does not necessarily result in a feasible schedule for a given sequence of pickup and drop-off locations. It is possible to reduce the user ride time by allowing delays in the service start time. Heuristic solution methods for the DARP usually apply the ``eight-step" method of \cite{cordeau2003tabu}, which constructs the feasible schedule by sequentially minimizing the possible violations of time windows, maximum route duration, and maximum user ride time. Recently, more advanced scheduling methods have been developed for the DARP, as in \cite{molenbruch2017multi} and \cite{bongiovanni2022ride}. 

As well as providing ride-sharing services, other recently trending approaches that help to reduce emissions and congestion include using Electric Vehicles (EVs) and developing autonomous driving technology. The employment of EVs offers the benefits of potentially fewer greenhouse gas emissions, lower energy cost per mile, and lower noise \citep{feng2013economic}. The introduction of autonomous driving leads to more flexibility in managing vehicle fleets, considerably lower operational costs, and better service quality \citep{fagnant2015operations,chen2016operations,burns2013transforming}. 
This article studies the Electric Autonomous DARP (E-ADARP), which was first introduced by \cite{bongiovanni2019electric}. Although the E-ADARP shares some of the constraints of the typical DARP (e.g., maximum user ride time, time window constraints), the E-ADARP is different from the typical DARP in two aspects:
(i) the employment of EAVs and a partial recharging policy, and
(ii) a weighted-sum objective that minimizes both total travel time and total excess user ride time;
The first aspect (i) requires checking battery feasibility for a given route, while the second aspect (ii) requires determining minimal-excess-time schedules for a feasible solution. 
The first aspect also implies other important features of the E-ADARP: (a) partial recharging is allowed en route, and (b) the maximum route duration constraints no longer exist due to the autonomy of vehicles. Allowing partial recharging introduces a trade-off between the time window and battery constraints: although longer recharging extends the driving range, it may also lead to time-window infeasibility for later nodes. Employing autonomous vehicles eliminates the need to predefine destination depots, as autonomous vehicles need to continuously relocate during their non-stop service. Other problem-specific constraints also increase the complexity of solving the E-ADARP. These constraints include a minimum battery level that must be maintained at the end of the route as well as limited visits to each recharging station. With these features and constraints, the possibility that a metaheuristic is trapped in local minima of poor quality increases, and feasible solutions are difficult to consistently find.

This paper offers a fourfold contribution. Firstly, 
we propose a new approach that efficiently computes minimum excess user ride time by introducing a fragment-based representation of paths.
Then, we apply an exact route evaluation scheme that executes feasibility checking in linear time. Combining these two methods, we propose an exact and efficient optimization of excess user ride time for an E-ADARP route. Secondly, we adapt a Deterministic Annealing (DA) algorithm to tackle the E-ADARP by integrating the proposed excess user ride time optimization method. To the best of our knowledge, this is the first time an exact excess user ride time optimization has been developed for computing locally optimal solutions within an algorithm for solving the E-ADARP. This method allows computing the minimum excess user ride time for a feasible E-ADARP route in linear time after preprocessing. 
Thirdly, we demonstrate the performance of the proposed DA algorithm through extensive numerical experiments. On the previously solved instances, the DA algorithm improves the solution quality by 0.16\% on average. We provide the best solutions for 70 out of 84 instances, among which 25 are new best solutions.
To further test our algorithm in solving large-scale instances, we construct new benchmark instances with up to 8 vehicles and 96 requests, and we provide 19 new solutions on newly-introduced instances. Finally, we extend the E-ADARP model to investigate the effects of allowing unlimited visits to recharging stations. The major difficulties for local search introduced by highly-constrained instances are lessened considering this more realistic situation. 


The remainder of this paper is organized as follows. Section \ref{sec::LiteratureRview} presents a comprehensive literature review on the DARP with Electric Vehicles (EVs) and Electric Vehicle Routing Problems (E-VRPs). Section \ref{sec::Problem} provides the problem definition and the notations of sets, parameters, and variables. It also discusses the objective function and constraints of the E-ADARP. Section \ref{sec::schedule} introduces the fragment-based representation of paths and the method to minimize total excess user ride time. A novel route evaluation scheme of linear time complexity is then described. Based on Section \ref{sec::schedule}, Section \ref{sec::DAalgo} presents the framework of the proposed DA algorithm and its main ingredients. In Section \ref{sec::result}, we conduct extensive computational experiments to demonstrate the performance of the proposed DA algorithm. This paper ends in Section \ref{sec::conclusion} with a summary of the results, contributions, and future extensions. 

\section{Literature Review}\label{sec::LiteratureRview}
The E-ADARP is a combination of the typical DARP and the E-VRP. However, it is distinct from these two contexts as it applies a weighted sum objective function that minimizes total travel time and total excess user ride time. This section briefly reviews the literature related to DARPs with EVs and E-VRPs. We emphasize works that apply heuristic and metaheuristic methods. We then review DARP-related articles that specifically focus on user ride time minimization. 

\subsection{Related literature of DARPs with EVs}
\cite{masmoudi2018dial} is the first work that introduces DARP with EVs. In their work, EVs are recharged through battery swapping and are assumed to have a constant recharging time. The authors use a realistic energy consumption model to formulate the problem and introduce three enhanced Evolutionary VNS (EVO-VNS) algorithm variants, which can solve instances with up to three vehicles and 18 requests. \cite{bongiovanni2019electric} considers EAVs in the DARP and introduces the E-ADARP. Partial recharging is allowed when vehicles visit recharging stations, and the authors impose a minimum battery level constraint for the vehicle's State of Charge (SoC) at the destination depot. The minimum battery level is formulated as $\gamma Q$, where $\gamma$ is the ratio of the minimum battery level to total battery capacity, and $Q$ is the total battery capacity. Three different $\gamma$ values are analyzed, i.e., $\gamma \in \{0.1,0.4,0.7\}$, meaning that 10\%, 40\%, and 70\% of the total battery capacity must be maintained at the destination depot. Solving the problem becomes more difficult when $\gamma$ increases. The authors formulate the problem into a three-index and a two-index model and introduce new valid inequalities in a Branch-and-Cut (B\&C) algorithm. When $\gamma = 0.1, 0.4$, the proposed B\&C algorithm obtains optimal solutions for 42 out of 56 instances. However, when $\gamma = 0.7$, the B\&C algorithm cannot solve 9 out of 28 instances feasibly, even with a two-hour run time. The largest instance that can be solved optimally by the B\&C algorithm contains 5 vehicles and 40 requests.
Recently, \cite{bongiovanni2022machine} have proposed a Machine Learning-based Large Neighborhood Search (MLNS) to solve the dynamic version of the E-ADARP. The proposed approach is a two-phase metaheuristic that sequential solves static E-ADARP subproblems. However, its performance on the previously defined static E-ADARP instances is not reported. Different from our algorithm, the authors focus on selecting destroy-repair operators at each iteration by a machine learning approach, which is trained offline on a large dataset produced through simulation. 

\subsection{Related literature of E-VRPs}
Extensive works have been conducted in the field of E-VRPs, e.g., \cite{erdougan2012green, schneider2014electric, goeke2015routing, hiermann2016electric, hiermann2019routing}. Among them, \cite{erdougan2012green} is the first to propose a Green VRP (G-VRP) using alternative fuel vehicles. These vehicles are allowed to visit a set of recharging stations during vehicle trips. The authors adapt two constructive heuristics to obtain feasible solutions and they further enhance these heuristics by applying local search. However, the proposed model does not consider capacity restrictions and time window constraints. \cite{schneider2014electric} propose a more comprehensive model named the Electric Vehicle Routing Problem with Time Windows (E-VRPTW). They extend the work of \cite{erdougan2012green} by using electric vehicles and considering limited vehicle capacity and specified customer time windows. They apply a Variable Neighborhood Search (VNS) algorithm hybridized by Tabu Search in local search to address E-VRPTW. The recharging stations are inserted or removed by a specific operator, and the recharged energy is assumed to be linear with the recharging time. They apply a full recharging policy on each visit to a recharging station. All the vehicles are assumed to be identical in terms of vehicle and battery capacity. \cite{goeke2015routing} extend the homogeneous E-VRPTW by considering a mixed fleet of electric and conventional vehicles. A realistic energy consumption model that integrates speed, load, and road gradient is employed. To address the problem, they propose an ALNS algorithm using a surrogate function to evaluate violations efficiently. 
\cite{hiermann2016electric} extend the work of \cite{goeke2015routing} by taking into account the heterogeneous aspect (i.e., fleet composition). They solve the problem by ALNS and determine the positions of recharging stations via a labeling algorithm. The recharging policy considered is also full recharging with a constant recharging rate.
\cite{hiermann2019routing} extend their previous study by considering partial recharging for a mixed fleet of conventional, plug-in hybrid, and electric vehicles. The engine mode selection for plug-in hybrid vehicles is considered as a decision variable in their study. A layered optimization algorithm is presented. This algorithm combines labeling techniques and a greedy route evaluation policy to calculate the amount of energy required to be charged and determine the engine mode and energy types. This algorithm is finally hybridized with a set partitioning problem to generate better solutions from obtained routes. More recently, \cite{lam2022branch} investigate a more practical case of E-VRPTW in which the availability of chargers at the recharging stations are considered. They propose a B\&C\&P algorithm which is capable of solving instances with up to 100 customers.

\subsection{Minimizing total or excess user ride time in DARPs }
\label{literature about minimizing excess user ride time}
There are several examples where a service-quality oriented objective is considered in the context of DARP (e.g., \cite{parragh2009heuristic,parragh2011introducing,molenbruch2017multi,bongiovanni2022ride}).  Among them, only three articles consider total user ride time/total excess user ride time as an objective. 
In the work of \cite{parragh2009heuristic}, a two-phase heuristic method is developed. A set of non-dominated solutions is constructed, minimizing a weighted sum of total distance traveled and mean user ride time under different weight combinations. In the route evaluation, the authors point out that the ``eight-step" method of \cite{cordeau2003tabu} does not aim to minimize the total user ride time. An increase in user ride time may happen when delaying the service start time at destination nodes. Therefore, they improve the original scheme of the ``eight-step" method by adapting the computation of forward time slack to avoid any increase in excess user ride time for requests served on a route. The resulting scheme is more restrictive in terms of feasibility and may lead to incorrect infeasibility declaration. 
This drawback is tackled in the scheduling heuristic proposed by \cite{molenbruch2017multi}. The heuristic starts by constructing a schedule (which may be infeasible) by setting the excess ride time of each request to its lower bound. Then, it gradually removes the infeasibility by shifting the service start time at some nodes while minimizing excess user ride time. However, the developed scheduling procedures in \cite{parragh2009heuristic} and \cite{molenbruch2017multi} are not proven optimal to minimize user ride time for a given route. \cite{bongiovanni2022ride} first proposes an exact scheduling procedure that can minimize the excess user ride time for a path without charging stations in polynomial time. 
Then, the authors extend the proposed scheduling procedure in the E-ADARP by integrating a battery management heuristic. However, the obtained schedules for an E-ADARP route are no longer exact as the excess-time optimal schedules may not be battery-feasible. The reported results show that on the investigated instances, the proposed scheduling procedure does not produce incorrect infeasible declarations, while others (i.e., \cite{cordeau2003tabu}, \cite{parragh2009heuristic}) do. 
To the best of our knowledge, no work in the literature can handle excess user ride time minimization exactly in the E-ADARP.

\subsection{Conclusion and proposed solution methodology}

From our review, we conclude that the effect of electric vehicles on the DARP has rarely been investigated in the previous literature. \cite{bongiovanni2019electric} is the only work that conducts a comprehensive study to optimize the static version of the DARP with EVs. However, the proposed B\&C algorithm requires important run-times and has difficulties providing high-quality solutions when solving medium- to large-sized instances, which limits its application in practice. The above limitation of \cite{bongiovanni2019electric} motivates us to propose an efficient metaheuristic algorithm that can provide high-quality solutions for E-ADARP instances within reasonable computational time. The efficiency of a metaheuristic largely depends on its neighborhood search mechanisms,  which perform a large number of evaluations. In the case of the DARP, these are route evaluations and cost computations. 
These two tasks are more complicated in the E-ADARP than in the DARP, as we allow partial recharging and minimize total excess user ride time for a given route. Existing scheduling procedures only obtain the approximation of minimum excess user ride time, which may deteriorate the solution quality and mislead search direction. Moreover, these procedures are time-consuming when applied in a metaheuristic as they are usually of quadratic time complexity and may introduce numerous repeated computations. Lastly, the battery constraints and a partial recharging policy increase the complexity of route evaluation in the E-ADARP. 

To overcome these issues, we propose an exact method of linear time complexity to compute the cost and evaluate the feasibility of an E-ADARP route based on battery-restricted fragments in Section \ref{sec::schedule}. Repeated computations are avoided via fragment enumeration in the preprocessing phase (Section \ref{sec::preproc}). These methods pave the way for an efficient DA algorithm (see Section \ref{sec::DAalgo}) and yield high-quality solutions for all instances (see Section \ref{sec::result}).

\section{The E-ADARP Description} \label{sec::Problem}

In this section, we present the mathematical notations of the E-ADARP that are originally introduced in \cite{bongiovanni2019electric} and are used throughout the paper. Then, we present the objective function and constraints of the E-ADARP. The final part discusses the practical interests of extending the original problem to allow unlimited visits to recharging stations.

\subsection{Notation and problem statement}

 The problem is defined on a complete directed graph $G=(V,A)$, where $V$ represents the set of vertices and $A$ is the set of arcs, i.e., $A = \{(i,j):i,j \in V, i \neq j\}$. Set $V$ can be further partitioned into several subsets, i.e., $V= N \cup S \cup O \cup F$, where $N$ represents the set of all customers, $S$ is the set of recharging stations, $O$ and $F$ denote the set of origin depots and destination depots, respectively. The set of all pickup vertices is denoted as $P =\{1,\cdots,i,\cdots,n\}$ and the set of all drop-off vertices is denoted as $D =\{n+1,\cdots,n+i,\cdots,2n\}$. The union of $P$ and $D$ is $N$, i.e., $N = P \cup D$. Each customer request is a pair $(i,n+i)$ for $i \in P$ and the maximum ride time for users associated with request $i$ is assumed to be $m_i$. A time window is defined on each node $i\in V$, denoted as $[e_i,l_i]$, in which $e_i$ and $l_i$ represent the earliest and latest time at which the vehicle starts its service, respectively. A load $q_i$ and a service duration $s_i$ is also associated for each node $i \in V$. 
 For pickup node $i \in P$, $q_i$ is positive. For the corresponding drop-off node $n+i$, we have $q_{n+i} = -q_i$. For other nodes $j \in O \cup F \cup S$, $q_j$ and $s_j$ are equal to zero. In this article, we tackle the static E-ADARP (i.e., all the customer requests are known at the beginning of the planning horizon $T_p$).

Each vehicle $k \in K$ must start with an origin depot $o \in O$ and end with a destination depot $f \in F$. In this study, the number of origin depots is equal to the number of vehicles, i.e., $|O| = |K|$, as in \cite{bongiovanni2019electric}. However, the set of destination depots can be larger than the set of origin depots, namely, $|F| \geqslant |O|$, which means a vehicle can select a depot from $F$ at the end of the route. An E-ADARP \emph{route} is defined as a path in graph $G$ originating from the origin depot and terminating in the destination depot, and passing through the pickup, drop-off, and charging station (if required) locations, which satisfies pairing and precedence, load, battery, time window, and maximum user ride time constraints. The E-ADARP consists in designing  $K$ routes, one for each vehicle, so that all customer nodes are visited exactly once, each recharging station and destination depot is visited at most once, and the weighted-sum objective function (presented in Section \ref{obj E-ADARP}) is minimized.
For unemployed vehicles, they travel directly from their designated origin depot to a destination depot. Vehicles are assumed to be heterogeneous in terms of their maximum vehicle capacities $C_k$ and homogeneous in terms of battery capacities (denoted as $Q$). 

The travel time on each arc $(i,j) \in A$ is denoted as $t_{i,j}$ and the battery consumption is denoted as $b_{i,j}$. We assume that $b_{i,j}$ is proportional to $t_{i,j}$ and we have $b_{i,j} = \beta t_{i,j}$, with $\beta$ being the energy discharging rate. When a vehicle recharges at a recharging station, the energy recharged is proportional to the time spent at the facilities. The recharging rate is denoted as $\alpha$. Energy units are converted to time units by defining $h_{i,j} = b_{i,j}/\alpha$. Then, the battery consumption $b_{i,j}$ on arc $(i,j)$ is converted to the time needed for recharging this amount of energy. 
Similarly, we can also convert the current energy level to the time needed to recharge to this energy level. Let $H$ denote the time required to recharge from zero to full battery capacity $Q$. Partial recharging is allowed while a vehicle visits recharging stations, and a minimum battery level $\gamma Q$ must be respected at destination depots. The triangle inequality is assumed to hold for travel times and battery consumption.

\subsection{Objective function of the E-ADARP } \label{obj E-ADARP}
A weighted sum objective is considered in this paper, which includes the total travel time for all the vehicles $k \in K$ and the total excess user ride time for all the customer requests $i \in P$. Equation (\ref{objective}) presents the formulation for the objective function. Considering the total excess user ride time in the objective function is also interesting, as it may help to improve the service quality by minimizing the total excess user ride time with no increase in the first objective if we consider the minimization in a strict lexicographical way. The objective function is:

\begin{equation} \label{objective}
     \min w_1\sum\limits_{k \in K}\sum\limits_{i,j \in V}t_{i,j}x_{i,j}^k + w_2\sum\limits_{i \in P}R_i
\end{equation}

where $x_{i,j}^k$ is a binary decision variable which denotes whether vehicle $k$ travels from node $i$ to $j$. $R_i$ denotes the excess user ride time of request $i \in P$ and is formulated as the difference between the actual ride time and direct travel time from $i$ to $n+i$. $w_1$ and $w_2$ are the weight factors for these two objectives.
We report in Table \ref{tab::notation} the notations and definitions for sets and parameters.

\begin{table}[ht]
\renewcommand\arraystretch{0.55}
    \caption{The E-ADARP problem sets, parameters notations and descriptions}
    \label{tab::notation}
    \begin{center}
    \begin{tabular}{c c}
    \toprule
    Sets& Definitions\\
    \hline
    $N = \{1,\cdots,n,n+1,\cdots,2n\}$& Set of pickup and drop-off nodes \\
    $P=\{1,\cdots,i,\cdots,n\}$& Set of pickup nodes \\
    $D=\{n+1,\cdots,n+i,\cdots,2n\}$& Set of drop-off nodes \\
    $K=\{1,\cdots,k\}$& Set of available vehicles \\
    $O=\{o_1,o_2,\cdots,o_k\}$& Set of origin depots \\
    $F=\{f_1,f_2,\cdots,f_h\}$& Set of all available destination depots (supposing the total number is $h$)\\
    $S=\{s_1,s_2,\cdots,s_g\}$& Set of recharging stations (supposing the total number is $g$) \\
    $V= N \cup S \cup O \cup F$& Set of all nodes \\
    \midrule
    Parameters& Definitions\\
    \hline
    $t_{i,j}$& Travel time from location $i\in V$ to location $j\in V$ \\
    $b_{i,j}$& Battery consumption from location $i\in V$ to location $j\in V$ \\
    $h_{i,j}$& The time needed for recharging $b_{i,j}, i,j\in V$ \\
    $e_i$& Earliest time at which service can begin at $i\in V$ \\
    $l_i$& Latest time at which service can begin at $i\in V$ \\
    $s_i$& Service duration at $i\in V$ \\
    $q_i$& Change in load at $i\in N$ \\
    $m_i$& Maximum user ride time for request $i\in P$ \\
    $C_k$& The vehicle capacity of vehicle $k$  \\
    $Q$& The battery capacity  \\
    $H$& Recharging time required to recharge from zero to $Q$ \\
    $\alpha$& The recharged energy per time unit\\
    $\beta$& The discharged energy per time unit\\
    $T_p$& Planning horizon\\
    $\gamma$& The ratio of minimum battery level at destination depot to $Q$ \\
    $w_1, w_2$ &Weight factors for total travel time and total excess user ride time \\
    \bottomrule
    \end{tabular}
    \end{center}
    \vspace{-4mm}
\end{table}

\subsection{Constraints of the E-ADARP } \label{constraints E-ADARP}

The E-ADARP consists of the following features that are different from the typical DARPs: 
\begin{itemize}
 \item [1)]
 Battery limitation and minimum battery level restriction, which introduce the detour to recharging stations;
 \item [2)]
 We allow partial recharging at recharging stations, and the recharging time must be determined;
 \item [3)]
 Vehicles locate at different origin depots and select the destination depot from a set of destination depots;
 \item [4)]
 Maximum route duration constraints are removed due to the autonomy of vehicles.
\end{itemize}

A solution of the E-ADARP is a set of $|K|$ routes and is called ``feasible" when all the following constraints are satisfied:
\begin{enumerate}
\item Every route starts from an origin depot and ends at a destination depot;
\item For each request, its corresponding pickup, and drop-off node  belong to the same route, and the pickup node is visited before its drop-off node;
\item User nodes and origin depots are visited exactly once, while each destination depot is visited at most once;
\item The maximum vehicle capacity must be respected at each node;
\item Each node is visited within its time window $[e_i,l_i]$ where $i \in V$. Vehicle can arrive earlier than $e_i$ but cannot arrive later than $l_i$. In the first case, waiting time occurs at $i$.
\item The maximum user ride time is not exceeded for any of the users;
\item The battery level at the destination depot must be at least equal to the minimal battery level;
\item The battery levels at any nodes of a route can not exceed the battery capacity and cannot be negative;
\item The recharging station can only be visited when there is no passenger on board;
\item Each recharging station $s \in S$ can only be visited at most once by all vehicles.
\end{enumerate}

An E-ADARP route is called ``feasible" if the above constraints, except for constraints (3) and (10), are fulfilled. Also, it should be noted that one can allow multiple visits to a recharging station by replicating the set of recharging stations, as in \cite{bongiovanni2019electric}.



\subsection{Multiple visits at recharging stations?}\label{sec::multiple?}

Each E-ADARP instance of \cite{bongiovanni2019electric} only contains a few recharging stations. In \cite{bongiovanni2019electric}, they first restrict the visit to the recharging station to at most one visit. Then, they investigate the effect of allowing multiple visits to recharging stations by replicating set $S$. Therefore, the number of visits to a recharging station must be predefined in their case. 
In our work, we relax this constraint and allow unlimited visits to the recharging stations in Section \ref{multiple section}.  As we have time window constraints and a minimal energy restriction at destination depots, visiting recharging stations more frequently leads to higher solution costs and increases the risk of violating time window constraints on succeeding nodes. 
We also conduct a sensitivity analysis on the maximum number of charging visits per station (denoted as $n_{as}$), and we perform our DA algorithm under different values of $n_{as}$ ($n_{as} = \{1,2,3,\infty\}$). 

\section{Excess User Ride Time Optimization} \label{sec::schedule}
The idea of our excess user ride time optimization method is as follows.
We first introduce a fragment-based representation of paths, which extends the one proposed in \cite{rist2021new} by additionally considering battery constraints for ensuring overall route feasibility in terms of energy consumption. Based on this representation of paths, each E-ADARP route can be represented by a series of battery-restricted fragments (see Definition \ref{battery restricted fragment}). 
 Then, we prove in Theorem \ref{theorem1} that the minimum total excess user ride time for a feasible route can be determined by summing the minimum excess user ride time of each \BRFrag. Following this idea, we enumerate all the feasible \BRFrags and calculate their minimum excess user ride times in the preprocessing phase (shown in Section \ref{sec::preproc}). With all the feasible fragments obtained as well as their minimum excess user ride time, we only need to check the feasibility of the route, which is realized via an exact route evaluation scheme of linear time complexity. 


\subsection{Representation of paths} \label{reprentation}
The most important characteristic of the E-ADARP is the incorporation of total excess user ride time in the objective function as well as the maximum user ride time in the constraints. Usually, the maximum user ride time constraints can be tackled by calculating forward time slack and delaying the service start time at some nodes (e.g., \cite{cordeau2003tabu}, \cite{kirchler2013granular}, \cite{parragh2009heuristic}). To minimize the total excess user ride time,  we declare one important point: total excess user ride time can only be minimized when vehicles finish their delivery (i.e., no open request on the path). We then introduce battery-restricted fragments: 

\begin{definition}[Battery-restricted fragment] \label{battery restricted fragment}
Assume that $\mathcal{F} = (i_1,i_2, \cdots,i_k)$ is a sequence of pickup and drop-off nodes, where the vehicle arrives empty at $i_1$ and leaves empty at $i_k$ and has passenger(s) on board at other nodes. Then, we call $\mathcal{F}$ a \BRFrag if there exists a feasible route of the form:
\begin{equation} \label{fragment definition}
    (o,s_{i_1},\cdots,s_{i_v},\overbrace{i_1,i_2, \cdots,i_k}^{\mathcal{F}},s_{i_{v+1}},\cdots,s_{i_m},f)
\end{equation}

where $s_{i_1},\cdots,s_{i_v},s_{i_{v+1}},\cdots,s_{i_m} (v,m \geqslant 0)$ are recharging stations, $o \in O$, and $f \in F$.
\end{definition}

\begin{figure}[!htp]
\centering
\includegraphics[width=16cm]{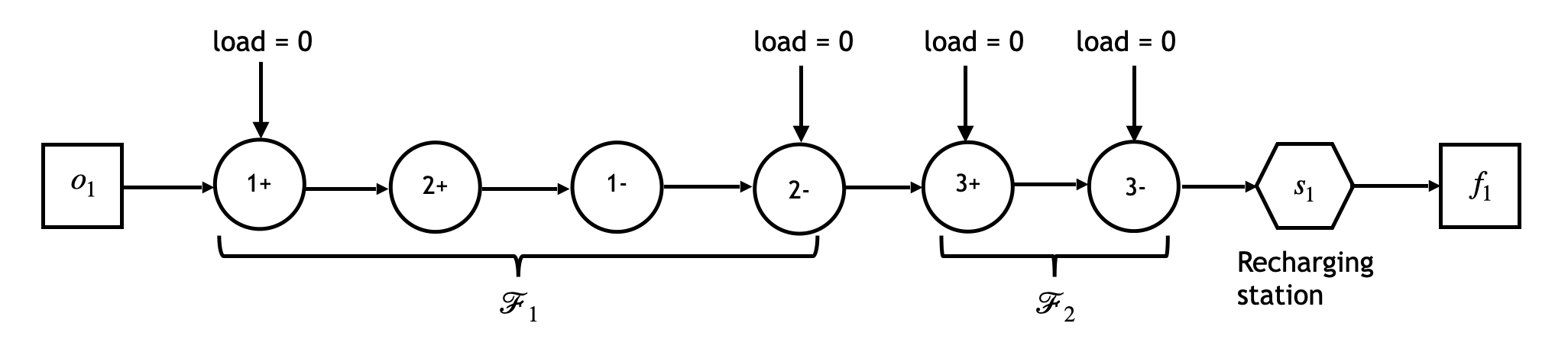}
\caption{\centering Example of \BRFrags}
\label{battery-restricted fragment example}
\end{figure}

It should be noted that if no recharging station is required in the route of Definition \ref{battery restricted fragment}, i.e., $v=m=0$ in Equation (\ref{fragment definition}), the battery-restricted fragment is equivalent to a fragment defined in \cite{rist2021new}. Figure \ref{battery-restricted fragment example} presents an example of a feasible route that consists of two battery-restricted fragments, i.e., $\mathcal{F}_1 = \{1+,2+,1-,2-\}$ and $\mathcal{F}_2 = \{3+,3-\}$. Note that $\mathcal{F}_1 \cup \mathcal{F}_2$ is not a battery-restricted fragment, as the vehicle becomes empty at intermediate node 2-. Based on this definition, each E-ADARP route can be regarded as the concatenation of several battery-restricted fragments, recharging stations (if required), an origin depot, and a destination depot.



Clearly, on each \BRFrag (hereinafter referred to as ``fragment"), the minimum excess user ride time can be calculated exactly.  We prove in the next section (Theorem \ref{theorem1}) that the minimum excess user ride time of route $\mathcal{R}$ can be calculated by summing the minimum excess user ride time on each fragment $\mathcal{F}_i \subseteq \mathcal{R}$. Then, we only focus on optimizing excess user ride time for each fragment.

\subsection{Excess user ride time optimization for a fragment} \label{ERT optimization}
Let $EU_{min}(\mathcal{R})$ and $EU_{min}(\mathcal{F})$ be the minimum excess user ride over route $\mathcal{R}$ and fragment $\mathcal{F}$, respectively. We have the following Theorem.
\begin{theorem} \label{theorem1}
If $\mathcal{R}$ is a feasible route and   $\mathcal{F}_1, \mathcal{F}_2, \cdots, \mathcal{F}_n$ are all the fragments on $\mathcal{R}$, then we have $EU_{min}(\mathcal{R}) = EU_{min}(\mathcal{F}_1) + EU_{min}(\mathcal{F}_2) + \cdots + EU_{min}(\mathcal{F}_n)$
\end{theorem}
We present the proof of Theorem \ref{theorem1} in \ref{proof of theorem}. Based on Theorem \ref{theorem1}, we convert the optimization of total excess user ride time for route $\mathcal{R}$ to the optimization of excess user ride time on its fragments $\mathcal{F} \subseteq \mathcal{R}$. Clearly, we can calculate the minimum excess user ride time directly if no waiting time is generated on a fragment. In the case of waiting time generated, one can compute the minimum excess user ride time if fragment $\mathcal{F}$ only contains a direct trip from one pickup node to the corresponding drop-off node. In the case that $\mathcal{F}$ contains two or more requests and waiting time generates for some $i \in \mathcal{F}$, the minimization of excess user ride time for $\mathcal{F}$ is equivalent to minimize a weighted sum of waiting time along $\mathcal{F}$, where weight factors are vehicle load at nodes with waiting time. To obtain the minimum excess user ride time, we resort to solving a Linear Program (LP), as presented in \ref{proof of theorem}.





Note that we ensure the maximum user ride time and vehicle capacity constraints when we generate fragments (will be explained in Section \ref{sec::preproc}). If a route $\mathcal{R}$ contains an infeasible fragment, it is discarded directly without further evaluation.

\subsection{Exact route evaluation scheme of linear time complexity} \label{exact route evaluation}
One challenge of the E-ADARP is tackling the trade-off between recharging time and time window constraints. A longer recharging time will extend the driving range and is beneficial to meet the energy restriction at the destination depot. However, the vehicle risks violating the time window constraints for the succeeding nodes. These two aspects interact, and it is hard to check the feasibility of a generated route (denoted as $\mathcal{R}$). We construct an exact route evaluation scheme of linear time complexity based on the forward labeling algorithm of \cite{desaulniers2016exact}. To the best of our knowledge, it is the first time an exact route evaluation scheme is developed to handle the DARP with EVs.

Given a route $\mathcal{R}$, we associate each node $i \in \mathcal{R}$ with a label $L_i :=\{(T_i^{rch_s})_{s \in S},T_i^{tMin},T_i^{tMax},T_i^{rtMax}\}$ including four resource attributes. We denote $\mathcal{P}_i$ as the partial path from the first node of $\mathcal{R}$ until node $i$. The definition of each resource attribute is shown as follows:
\begin{enumerate}
\item $T_i^{rch_s}$: The number of times recharging station $s \in S$ is visited along  $\mathcal{P}_i$;
\item $T_i^{tMin}$: The earliest service start time at vertex $i$ assuming that, if a recharging station is visited prior to $i$ along  $\mathcal{P}_i$, a minimum recharge (ensuring the battery feasibility up to $i$) is performed;
\item $T_i^{tMax}$: The earliest service start time at vertex $i$ assuming that, if a recharging station is visited prior to $i$ along  $\mathcal{P}_i$, a maximum recharge (ensuring the time-window feasibility up to $i$) is performed;
\item $T_i^{rtMax}$: 
The maximum recharging time required to fully recharge at vertex $i$ assuming that, if a recharging station is visited prior to $i$ along $\mathcal{P}_i$, a minimum recharge (ensuring the battery feasibility up to $i$) is performed;
\end{enumerate}

The initial label is defined as $\{ (\overbrace{0,\cdots,0}^{|S| \text{ times}}),0,0,0\}$. We compute the succeeding label $L_j$ from the previous label $L_i$ by Resource Extension Functions (REFs): 

\begin{equation}
T_j^{rch_s}= T_i^{rch_s} +
\begin{cases}
    1,\quad & \text{if $j = s$} \\
    0, \quad & \text{otherwise}
\end{cases}
\end{equation}

\begin{equation}
T_j^{tMin}=
\begin{cases}
   \max \{e_j,T_i^{tMin}+t_{i,j}+s_i\} , \quad &\text{if $T_i^{rch}= \emptyset$}  \\
   \max \{e_j,T_i^{tMin}+t_{i,j}+s_i\}+Z_{i,j}, \quad &\text{otherwise} 
\end{cases}
\end{equation}

\begin{equation}
T_j^{tMax}=
\begin{cases}
  \min \{l_j,\max\{e_j,T_i^{tMin}+T_i^{rtMax}+t_{i,j}+s_i\}\}, \quad & \text{if $i \in S$} \\
   \min \{l_j,\max\{e_j,T_i^{tMax}+t_{i,j}+s_i\}\}, \quad & \text{otherwise} 
\end{cases}
\end{equation}

\begin{equation}
T_j^{rtMax}=
\begin{cases}
  T_i^{rtMax}+h_{i,j}, \quad & \text{if $T_i^{rch} = \emptyset$} \\
  \min\{H,\max\{0,T_i^{rtMax}-S_{i,j}\}+h_{i,j}\}, \quad & \text{otherwise}
\end{cases}
\end{equation}

where:

\begin{equation}
S_{i,j}(T_i^{tMin},T_i^{tMax},T_i^{rtMax})=
\begin{cases}
  \max\{0, \min\{e_j-T_i^{tMin}-t_{i,j}-s_i,T_i^{rtMax}\}\},\quad & \text{if $i\in S$ }\\
  \max\{0, \min\{e_j-T_i^{tMin}-t_{i,j}-s_i,T_i^{tMax}-T_i^{tMin}\}\},\quad & \text{otherwise} 
\end{cases}
\end{equation}

\begin{equation}
    Z_{i,j}(T_i^{tMin},T_i^{tMax},T_i^{rtMax}) = \max\{0, \max\{0, T_i^{rtMax}-S_{i,j}(T_i^{tMin},T_i^{tMax},T_i^{rtMax})\}+h_{i,j}-H\}
\end{equation}

The $S_{i,j}$ is the slack time between the earliest time window $e_j$ at $j$ and the earliest arrival time to $j$. If $i$ is a recharging station, $S_{i,j}$ is the maximum amount of recharging time that can be performed at $i$, namely $T_i^{tMax}-T_i^{tMin}$. $Z_{i,j}$ is the minimum recharging time required to keep battery feasibility accounting for the available slack at the previous recharging station. 

According to \cite{desaulniers2016exact}, we have following proposition:
\begin{proposition} \label{prop1}
The route $\mathcal{R}$ is feasible if and only of $ \forall j \in R $, the label $L_j$ satisfies:

\begin{equation} \label{tw at j}
    T_j^{tMin} \leqslant l_j,\quad T_j^{tMin} \leqslant T_j^{tMax},\quad  T_j^{rch_s} \leqslant 1,\quad T_j^{rtMax} \leqslant
    \begin{cases}
      (1-\gamma)H, \quad & j \in F\\
      H, \quad & \text{otherwise}\nonumber
    \end{cases}
\end{equation}

\end{proposition}

Clearly, the feasibility checking algorithm is of linear time complexity with respect to the length of the input route. After checking the feasibility, the total cost of route $\mathcal{R}$ is obtained by summing the travel time of arcs and the excess user ride time of fragments, recalling Theorem \ref{theorem1}.

\section{Deterministic Annealing Algorithm for the E-ADARP} \label{sec::DAalgo}

Based on Section \ref{ERT optimization} and Section \ref{exact route evaluation}, we establish a DA algorithm that ensures minimal excess user ride time for a generated solution and integrates an exact route evaluation. Different types of local search operators are embedded in the proposed DA algorithm to solve the E-ADARP. 


DA was first introduced by \cite{dueck1990threshold} as a variant of simulated annealing. Recent research shows that DA can obtain near-optimal or optimal solutions for a series of vehicle routing problems \citep{braysy2008effective,braekers2014exact}. To the best of our knowledge, the only paper that implements DA to solve the DARP is that of \cite{braekers2014exact}. Applying DA algorithm provides several advantages, and the most important one is its easy parameter tuning process, as the DA algorithm mainly relies on a single parameter. In addition, the DA algorithm is proved to be very efficient in solving the typical DARP. However, \cite{braekers2014exact} considers a single-objective case in the DARP. To solve the E-ADARP, we adapt the DA algorithm to accommodate problem-specific features of the E-ADARP by integrating the proposed excess user ride time optimization approach.

The framework for the proposed DA algorithm is depicted in Algorithm \ref{alg:meta-heuritic}.
The algorithm input is an initial solution $x_{init}$ constructed by a parallel insertion heuristic (presented in Section \ref{4.2}) and the initial settings of DA-related parameters. These parameters include: (i) a maximal number of iterations $N_{iter}$; (ii) the initial and maximal temperature $\Theta_{max}$; (iii) restart parameter $n_{imp}$. It should be mentioned that the initial solution $x_{init}$ is feasible for the E-ADARP constraints, except that only a subset of requests may be served. The solution cost of the initial solution is denoted as $c(x)$, and the number of requests served in the initial solution is updated to $N_{req}$ so that a lexicographic optimization considers cost comparison in $c(x)$ values only if it does not worsen the number of requests served. A list of indexed operators $opt_1, \dots, opt_z $ are operated sequentially in each DA iteration (presented in Section \ref{4.4}). 

\begin{algorithm}
  \caption{DA Algorithm for the E-ADARP}  
  \label{alg:meta-heuritic}
  \begin{algorithmic}[1]
  \small
    \Require
     Initial solution $x_{init}$, initial values of $N_{iter}$,  $\Theta_{max}$, and $n_{imp}$. $\Theta$ is set to $\Theta_{max}$
    \Ensure
     Best solution $x_b$ found by our algorithm;
    \While{$iter \leqslant N_{iter}$}
    \State $i_{imp} \leftarrow i_{imp} + 1$;
    \For{$j=1 \rightarrow z-1$}
    \State Apply local search operator $opt_j$ on $x$ to obtain neighborhood solution $x'$;
    \If{$c(x') < c(x) + \Theta $}
        \State $x \leftarrow x'$;
    \EndIf
    \EndFor
    \If{$N_{req} < n$}
    \State Apply $opt_z$ operator to add request to generate neighborhood solution $x'$;
    \State Update the number of requests served in $x'$ as $N_{req}'$;
    \EndIf
    \If{$(c(x')<c(x_b)$ \textbf{and} $N_{req}' = N_{req})$ \textbf{or} $N_{req}' > N_{req}$} 
    \State $x_b \leftarrow x'$
    \State $i_{imp} \leftarrow 0$
    \Else{}
     \State $\Theta \leftarrow \Theta- \Theta_{max}/\Theta_{red}$
     \If{$\Theta<0$}
     \State $r \leftarrow$ random number between 0 and 1
     \State $\Theta \leftarrow r \times \Theta_{max}$
    \If{$i_{imp} > n_{imp}$}
    \State $x \leftarrow x_b$
    \State $i_{imp} \leftarrow 0$
    \EndIf
     \EndIf
    \EndIf
    \State $iter \leftarrow iter + 1$
    \EndWhile
    \State \textbf{return $x_b$}
  \end{algorithmic}  
\end{algorithm}

There are two steps in the algorithm: local search and threshold update. At the beginning of the algorithm, the threshold value $\Theta$ is set to $\Theta_{max}$, and the best solution $x_b$ and current solution $x$ is initialized to an initial solution $x_{init}$. 
During the local search process, local search operators are applied to alter the current solution. In the next step, the threshold value is updated and restarted when the value is negative.

In the local search process, 
we first remove the existing recharging stations on the current route and then generate a random neighborhood solution $x'$ from the current solution $x$ by applying different operators.
In the case of neighborhood solution $x'$ satisfies $c(x') < c(x) +\Theta$ but violates battery constraints, 
we call an insertion algorithm to repair $x'$ by inserting recharging stations at proper places (presented in Section \ref{4.3}). Solution $x'$ is accepted to become the new current solution when the number of assigned requests increases or the total cost is less than that of the current solution plus the threshold value $\Theta$. 

In the threshold update process, when no new global best solution is found, $\Theta$ is reduced by $\Theta_{max}/\Theta_{red}$, where $\Theta_{red}$ is a predefined parameter. To ensure that $\Theta$ is always non-negative, we reset $\Theta$ to $r \times \Theta_{max}$, with $r$ a random number generated between zero and one whenever $\Theta$ becomes negative. The search is restarted from $x_b$ when no improvement is found in $n_{imp}$ iterations and $\Theta$ becomes negative.

\subsection{Parallel insertion heuristic} \label{4.2}
While in most of the literature, the initial solution is often generated randomly, we construct our initial solution by a parallel insertion algorithm considering the time window and spatial closeness, as in \cite{masmoudi2017hybrid}. First, we sort all the requests $(i,n+i), i \in P$ in increasing order based on $e_i$. Then, we randomly initialize $k$ routes $\{\mathcal{R}_1, \cdots, \mathcal{R}_k\}$ ($0 < k \leqslant K$ with $K$ being the number of total vehicles). Each of the $k$ first requests in the sorted request list are assigned randomly to different routes. These requests are deleted from the list of requests. 

Then, we sort the route list $\{\mathcal{R}_1, \cdots, \mathcal{R}_k\}$ in increasing order with regards to the distance between the last node of the analyzed route and the pickup node of the first request remaining in the request list. The first request is assigned to the first route in the route list. To insert the selected request, we enumerate all the possible insertion positions and insert the corresponding pickup node and drop-off node in a feasible way on this route. If this request cannot be inserted feasibly, then we move to the second route. This process is repeated until this request is inserted or all the routes are analyzed. If this request cannot be inserted in any of the existing routes, we move to the second request in the list and repeat the above process. After this process, if some requests are still not assigned, a new route is activated, and the above process will be repeated. The algorithm terminates when the request list is empty or the existing requests in the list cannot be inserted into any of the routes in a feasible way.

\subsection{Recharging station insertion for a given route} \label{4.3}
If a route $\mathcal{R} \in x'$ only violates the battery constraints and neighborhood solution $x'$ has $c(x') < c(x)+T$, we insert a/several recharging station(s) to repair $\mathcal{R}$. 
First, we create two empty sets, one is to store repaired route candidates (called ``list of feasible routes"), the other is to store potential route candidates (called ``list of candidate routes"). 
For each possible insertion position, we select a random recharging station from the set of available stations to insert. We do not consider inserting the best station (e.g, the closest one), as we may have other battery-infeasible routes in $\mathcal{R}$, which requires visiting this recharging station to be repaired.
If a feasible route is generated after insertion, we add it to the list of feasible routes. Otherwise, we store this route in the list of candidate routes. Suppose the route is still infeasible after trying all the possible insertion positions. In that case, we move to the next iteration to insert another recharging station for all the possible positions of all the candidate routes. The algorithm returns the repaired minimum-cost feasible route if $\mathcal{R} $ can be repaired or an empty set otherwise. 
For acceleration, we only consider repairing the route containing less than $N_{rch}$ recharging stations and we take $N_{rch} = \lceil |S|/2 \rceil$.

\subsection{Local search} \label{4.4}

We design seven operators (i.e., $opt_1, \cdots, opt_7$ in Algorithm \ref{alg:meta-heuritic}) to improve the initial solution generated from the constructive heuristic.
Among them, three are intra-route operators (i.e., \textit{ex-pickup}, \textit{ex-dropoff}, and \textit{ex-2-neighbor}), three are inter-route operators (i.e., \textit{2-opt}, \textit{relocate}, and \textit{exchange}). 
The last operator named \textit{add-request} is applied in each iteration on neighborhood solution $x'$, which is generated after applying $opt_1, \cdots, opt_5$, if there exists un-served requests.

\subsubsection{Intra-route operators}



\begin{figure}[H]
    \centering
    \subfigure[Ex-pickup operator example.]{
    \includegraphics[width=7cm]{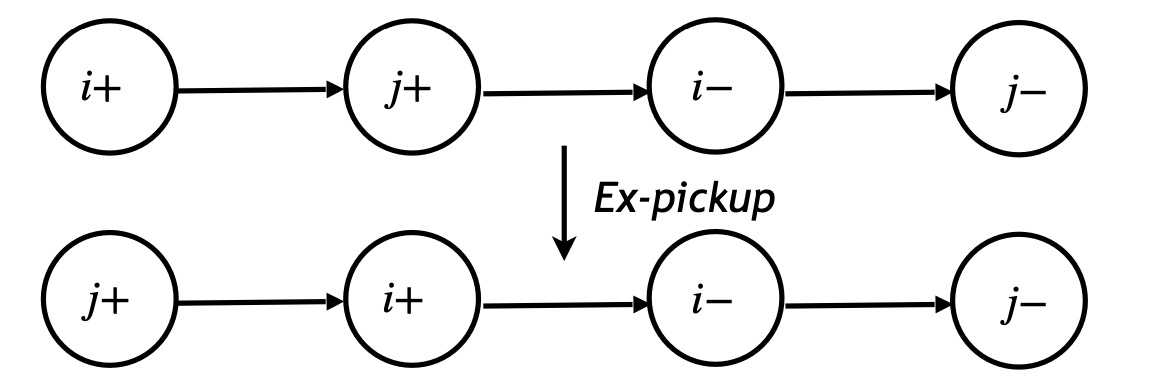}
    \label{exchange pickup}
    }
    \quad
    \subfigure[Ex-dropoff operator example.]{
    \includegraphics[width=6.5cm]{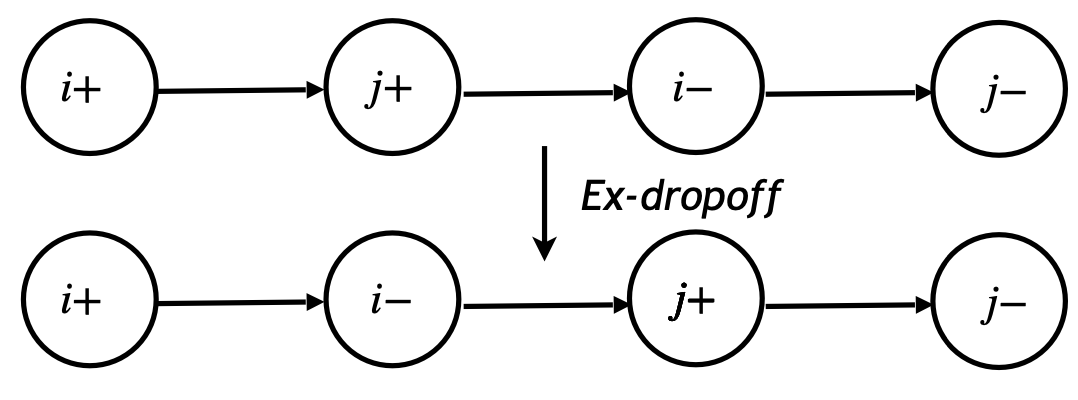}
    \label{exchange dropoff}
    }
    \quad
    \subfigure[Ex-2-neighbor operator example.]{
    \includegraphics[width=7cm]{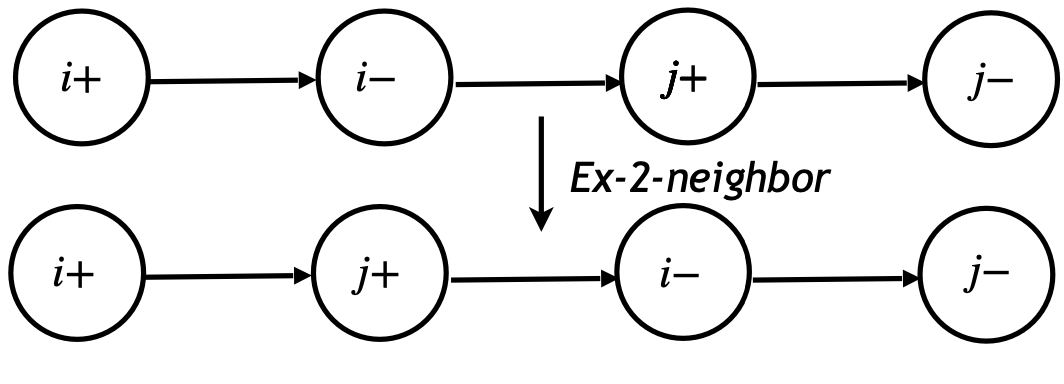}
    \label{exchange two neighbor}
    }
    \caption{\centering Intra-route operators example}
\end{figure}

\textit{Ex-pickup} operator swaps the positions of two consecutive nodes $(i+,j+)$, where node $i+$ is a pick-up node and node $j+$ is not the corresponding drop-off node. An example is shown in Figure \ref{exchange pickup}. In each iteration, one pick-up node is selected randomly. If the successor of this pick-up node does not correspond to its drop-off node, then the two positions are exchanged.

\textit{Ex-dropoff} operator creates a neighborhood solution by swapping the positions of two consecutive nodes $(j+,i-)$, where point $i-$ is a drop-off node and point $j+$ is not the corresponding pick-up node. Figure \ref{exchange dropoff} shows an example of how \textit{ex-dropoff} works. In each iteration, one drop-off node is selected randomly, if the precedent node of this drop-off node does not correspond to its pick-up node, then the two positions are exchanged.

There is another situation shown in Figure \ref{exchange two neighbor}, where the successor of pick-up node $i+$ is its drop-off $i-$, and the predecessor of drop-off node $j-$ is its corresponding pick-up $j+$, but we can still exchange $i$- and $j+$ to create a new neighborhood solution. This operation is realized by \textit{ex-2-neighbor} operator.


\subsubsection{Inter-route operators}

Inspired by state-of-the-art heuristic methods, we apply three widely-used inter-route operators to generate neighbors of the current solution, as in \cite{braekers2014exact}. 
\textit{Two-opt} operator selects two random routes and splits each route into two parts by a randomly selected zero-split node $i$ such that $i \in D \cup S$. Then, the first part of the first route is connected with the second part of the second route and the first part of the second route is connected with the second part of the first route. Note that \textit{2-opt} is able to realize the exchange of several requests at one iteration. 
\textit{Relocate} operator randomly removes one request from a random route and re-inserts the request at the best position of another route. The best position means the position that brings the least increase on solution cost after inserting the selected request. 
\textit{Exchange} operator 
swaps two random requests of two randomly-selected routes. The selected requests are re-inserted into the best position of the other route. 


\subsubsection{Insertion operator}

\begin{figure}[H]
\centering
\includegraphics[width=10cm]{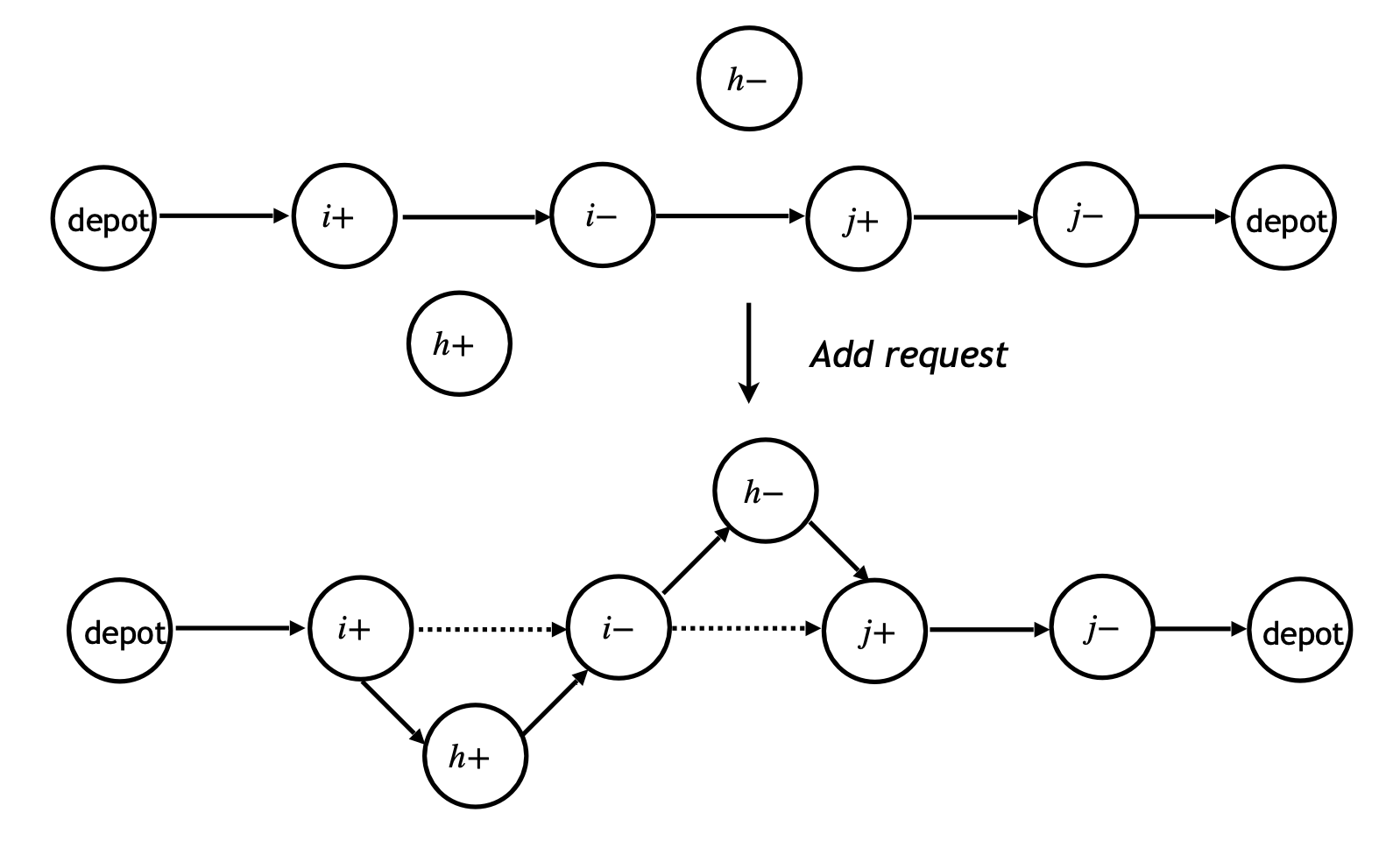}
\caption{\centering \textit{Add-request} operator example}
\label{addNewRequest}
\end{figure}

\textit{Add-request} operator is applied in each iteration when there exist uninserted requests for current solution $x$. 
This operator tries to insert one uninserted request into a random route of $x$. When all the requests are served in $x$, this operator will no longer be applied. Figure \ref{addNewRequest} describes how \textit{add-request} adds uninserted request $(h+,h-)$ on a route.

\subsection{Implementation details} \label{sec::preproc}
This section presents the preprocessing works and the algorithm implementation details for allowing multiple/unlimited visits to recharging stations. The preprocessing works include: time window tightening, arc elimination, and fragment enumeration. 

\subsubsection{Preprocessing works}
We first introduce two traditional methods introduced by \cite{cordeau2006branch}, which includes time window tightening and arc elimination. Then, we introduce fragments enumeration method.


Time window tightening is executed as:

\begin{itemize}
 \item For $i \in P$, $e_i$ is set to $\max\{ e_i,e_{n+i}-m_i-s_i \}$ and $l_i=\min \{ l_{n+i}-t_{i,n+i}-s_i,l_i \}$; 
 \item For $i \in D$, $e_{n+i}= \max\{e_{n+i},e_i+t_{i,n+i}+s_i \}$, and $l_{n+i}=\min\{l_i+m_i+s_i,l_{n+i}\}$.
 \item For $s \in S$, the time window can be tightened by considering the travel time from the origin depot to recharging station and from recharging station to the destination depot. The earliest time to start service at charging station $s$ is set to $\min\{e_j+t_{j,s} \}$, $\forall j \in O$; the latest time at charging station $s$ to start service at recharging station is $\max\{T_p-t_{s,j} \}, \forall j \in F$;
  \item For $i \in O \cup F$, the earliest time window $e_i$ is set to $\max\{0,\min\{e_j-t_{i,j} \}\}, \forall j \in P$, and $l_i=\min\{l_i,\max\{l_j+s_i+t_{j,i} \}\}, \forall j \in D$.
\end{itemize}

 It should be noted that we tighten the value of $e_i$ for node $i \in P$ by considering the earliest possible service start time at node $i \in P$ for arriving at the corresponding drop-off node $n+i$ at $e_{n+i}$. Hence, $e_i = \max\{ e_i,e_{n+i}-m_i-s_i \}$. Similarly, we tighten $l_i, i \in P$ by considering the latest possible service start time at node $i \in P$ for arriving at node $n+i$ at $l_{n+i}$, i.e., $l_i=\min \{ l_{n+i}-t_{i,n+i}-s_i,l_i \}$.
 The arc elimination process follows the method of \cite{cordeau2006branch}. We reduce the number of arcs in the graph by removing arcs that will not lead to a feasible solution. 

We further accelerate computations by enumerating all feasible fragments before computation, as in \cite{alyasiry2019exact}, \cite{rist2021new}. This method simplifies route evaluation and avoids recalculations as we only need to query information from each fragment. We enumerate all the feasible fragments with depth-first search and calculate their minimum excess user ride time. Then, the total excess user ride time of a route $\mathcal{R}$ can be calculated by summing $EU_{min}(\mathcal{F}), \mathcal{F} \subseteq \mathcal{R}$, recalling Theorem \ref{theorem1}.

To generate all feasible fragments, we start from each pickup node and extend it node by node respecting time window, capacity, battery, and maximum user ride time constraints. We assume that the vehicle starts from each pickup node with a full battery level. The maximum user ride time, vehicle capacity constraints are checked during the extension process. For each node on a fragment, it must have a positive battery level.

Note that if a fragment only contains one request, we calculate the excess user ride time directly and check the maximum user ride time constraints. If a fragment contains two or more requests, we resort to a LP solver (Gurobi) to solve the LP model (shown in \ref{proof of theorem}) and check the maximum user ride time constraints. For each feasible fragment, the obtained minimum excess user ride time value is recorded.
In \ref{preliminary for seg enumeration}, we conduct a preliminary test and provide details for fragment enumeration on each instance. For all the instances, the fragment enumeration can be fulfilled in a matter of seconds. In the computational experiments, we report the CPU time which includes the computational time for performing all the preprocessing works in Section \ref{sec::result}.

\subsubsection{Adapt DA algorithm to allow multiple visits to each recharging station}
\label{adapt DA to allow multiple visits}
Different from the model of \cite{bongiovanni2019electric}, we allow multiple visits at each recharging station without the need to replicate set $S$ in Section \ref{multiple section}.
In the case of $n_{as} = 2,3$, we replicate the recharging station set $S$ to allow at most two and most three visits per station. All the ingredients remain the same in these two cases. In the case of $n_{as} = \infty$, we remove the feasibility checking rule $T_j^{rch_s} \leqslant 1$ to allow one route visiting multiple times for a station. When selecting a recharging station to insert in a route, we relax the set of available recharging stations to $S$. This operation allows inserting a recharging station that has already been used in other routes.

\section{Computational Experiments and Results} \label{sec::result}
In this section, we conduct extensive numerical experiments and analyze the results. All algorithms are coded in Julia 1.7.2 and are performed on a standard PC with an Intel Xeon Gold 6230 20C at 2.1GHz. 
This section is organized as follows. The benchmark instances for the computational experiments and abbreviations used in the Tables are introduced in the first part. Then, a sensitivity analysis is conducted to find good parameter settings for the proposed DA algorithm in Section \ref{DA parameters}. After ensuring the robustness of parameters and operators, we validate the performance of the proposed algorithm on the standard E-ADARP instances compared to the state-of-the-art results in Section \ref{DA performance}. Section \ref{multiple section} investigates the effect of allowing multiple visits to recharging stations.

\subsection{Benchmark instances and abbreviations}

This section presents the benchmark instances used to test the algorithm performance, their characteristics, and the notations for the computational experiments.

\subsubsection{Benchmark Instances}
Instances are named following the pattern xK-n-$\gamma$, where $K$ is the number of vehicles, $n$ is the number of requests, and $\gamma \in \{0.1,0.4,0.7\}$.
Three sets of instances are considered in the experiments, which differentiate by $x \in \{a,u,r\}$: 
\begin{itemize}
 
\item “a” denotes the standard DARP benchmark instance set  from \cite{cordeau2006branch} extended with features of electric vehicles and recharging stations by \cite{bongiovanni2019electric}. To simplify, we call them type-a instances. For type-a instances, the number of vehicles is in the range $2\leqslant K \leqslant 5$, and
the number of requests is in the range $16\leqslant n \leqslant 50$.
\item “u” denotes instances based on the ride-sharing data from Uber Technologies (instance name starts with “u”) that were adopted from \cite{bongiovanni2019electric}. To simplify, we call them type-u instances.
For type-u instances, the number of vehicles is in the range $2\leqslant K \leqslant 5$, and
the number of requests is in the range $16\leqslant n \leqslant 50$, as in type-a instances.
\item “r” denotes larger DARP benchmark instances built from \cite{ropke2007models} using the same extension rules
to have E-ADARP instances from DARP instances. To simplify, we call them type-r instances.
For type-r instances, the number of vehicles  is in the range $5\leqslant K \leqslant 8$ and
the number of requests is in the range $60\leqslant n \leqslant 96$.

\end{itemize}

Type-a instances are supplemented with recharging station ID, vehicle capacity, battery capacity, the final state of charge requirement, recharging rates, and discharging rates. The same operation is applied to type-r instances to generate a large-scale set of instances. The vehicle capacity is set to three passengers, and the maximum user ride time is 30 minutes. As in \cite{bongiovanni2019electric}, recharging rates and discharging rates are all set to 0.055KWh per minute according to the design parameter of EAVs provided in: \url{https://www.hevs.ch/media/document/1/fiche-technique-navettes-autonomes.pdf}. The efficient battery capacity is set to 14.85 KWh, and the vehicle can approximately visit 20 nodes without recharging. 

The ride-sharing dataset of Uber is obtained from the link: \url{https://github.com/dima42/uber-gps-analysis/tree/master/gpsdata}. Type-u instances are created by extracting origin/destination locations from GPS logs in the city of San Francisco (CA, USA) and applying Dijkstra’s shortest path algorithm to calculate the travel time matrix with a constant speed setting (i.e., 35km/h). Recharging station positions can be obtained through Alternative Fueling Station Locator from Alternative Fuels Data Center (AFDC). For a more detailed description of instances development, the interested reader can refer to \cite{bongiovanni2019electric}. The preprocessed data that extract requests information from the raw data provided by Uber Technologies are published on the website (\url{https://luts.epfl.ch/wpcontent/uploads/2019/03/e_ADARP_archive.zip}).

Following \cite{bongiovanni2019electric}, we consider three different $\gamma$ values, i.e., $\gamma \in \{0.1,0.4,0.7\}$, representing different minimal battery restrictions at the destination depot. For weight factors, we take $w_1 = 0.75$ and $w_2 = 0.25$.

 \subsubsection{Abbreviations in the tables}

The DA algorithm has deterministic rules to accept a solution and the sequence of neighborhoods, which is contrary to Simulated Annealing. There remains a randomized part in the selection of neighboring solutions. Unless indicated,  we perform 50 runs on each instance with different seeds to analyze the statistical distribution of the solution quality.

For each instance, we present the following values:

\begin{itemize}
 \item $BC'$ is the cost of best solutions from B\&C algorithm reported in \cite{bongiovanni2019electric};
 \item $BC$ is the cost of best solutions found by the proposed DA algorithm over 50 runs;
 \item $AC$ is the average-cost solution found by the proposed DA algorithm over the 50 runs.
 \item $Q1$ is the middle number between the best-obtained solution and the median of all the solutions over 50 runs;
 \item $Q3$ is the middle number between the median of all the solutions over 50 runs and the worst solutions yielded;
\end{itemize}

To analyze the distribution of the solution found for the 50 runs, we calculate solutions gaps to $BC'$.
Assuming a solution with value $v$ ($v$ could be $BC$, $Q1$, $Q3$), we compute its gap to $BC'$ by:

\begin{equation}
    gap =  \dfrac { v - BC'} {BC'} \times 100\% \nonumber
\end{equation}

Note that type-r instances for the E-ADARP are studied here for the first time, we therefore replace $BC'$ with $BC$ in the above formula to analyze the gaps of $Q1$/$AC$/$Q3$ to $BC$. 

We present the following average values to analyze the consistency of the proposed DA algorithm:
 
\begin{itemize}

\item $Q1\%$ is the average  gap to $BC'$ of the first quartile value over the different runs;
\item $Q3\%$ is the average  gap to $BC'$ of the third quartile value over the different runs;
\item $BC\%$  is the average gap of $BC$ to $BC'$ over the different runs;
\item  $AC\%$ is the average gap of $AC$ to $BC'$ over the different runs;
\item FeasRatio is the ratio of feasible solutions found among all the solutions generated by DA algorithm;
\item CPU is the average computational time of the DA algorithm (preprocessing time is included) in seconds;
\item CPU$'$ is the computational time of the B\&C algorithm reported in \cite{bongiovanni2019electric} in seconds;
\item NC (Not Calculable) means that there are unsolved instances under the analyzed parameter and we cannot calculate gaps.
\item NA (Not Available) indicates that corresponding value (e.g., $BC$, $BC'$) is not available as the analyzed algorithm cannot provide a feasible solution. 
\item A dash ``--" indicates that the DA algorithm finds new best solutions on a previously unsolved instance and we cannot calculate the gap. 
\end{itemize}
In Section \ref{multiple section}, we present DA algorithm results when allowing multiple visits to each recharging station. To distinguish, subscripts ``2", ``3", and ``$\infty$" are added to $BC$, $AC$, and CPU to denote $n_{as} = 2,3, \infty$, respectively. As \cite{bongiovanni2019electric} provides results on type-u instances with $n_{as} = 2,3$, we add their reported results in the column named $BC_2'$ and $BC_3'$ of Table \ref{mul uber} and compare our DA algorithm results to theirs.

 
\subsection{Parameter tuning for the DA algorithm}
\label{DA parameters}

The performance of the proposed algorithm depends on several parameters that must be set in advance. To ensure the algorithm's performance, we first identify robust parameter settings. We analyze different settings of parameters on the type-a instance set, as it contains instances of different sizes and is enough to select good parameters. 
For a comprehensive overview, we take into account different scenarios, i.e., $\gamma = 0.1, 0.4, 0.7$, for each parameter setting. 

The DA-related parameters are:

\begin{itemize}
 \item Number of iterations $N_{iter}$ ;
\item  Maximum threshold value $\Theta_{max}$;
\item Threshold reduction value $\Theta_{red}$;
\item Restart parameter $n_{imp}$.
 \end{itemize}

To avoid re-tuning $\Theta_{max}$ when using different instances, we use a relative value for $\Theta_{max}$. The maximum threshold value is expressed as the product of the average distance between two nodes in the studied graph (denoted $\Bar{c}$) and a predefined parameter $\theta_{max}$, that is $\Theta_{max}= \Bar{c} \times \theta_{max}$, where $\theta_{max}$ is initially set to 1.5. For other parameters like $\Theta_{red}$ and $n_{imp}$, we take the same settings as in \cite{braekers2014exact}: $\Theta_{red}=300$ and $n_{imp}=50$. 

\subsubsection{Sensitivity analysis and parameter tuning for $\theta_{max}$}

The sensitivity analysis results for $\theta_{max}$ under $\gamma = 0.1, 0.4, 0.7$ are shown in \ref{detailed parameter tuning}, and we test seven values for $\theta_{max}$. For each value of $\theta_{max}$, we perform ten runs on each instance and iterate the proposed algorithm 10000 times for each run. Under each energy restriction, we report $BC\%$, $AC\%$, $Q1\%$, $Q3\%$ over ten runs for the analyzed $\theta_{max}$ value. For the scenario of $\gamma = 0.7$, we report FeasRatio and average CPU time. We present detailed results on each instance under different settings of $\theta_{max}$ in \ref{detailed parameter tuning}.

From Table \ref{t_max}, in the case of $\gamma = 0.1$, the algorithm performs well under all the settings of $\theta_{max}$. Among them, 0.6 seems to be the best with regard to gap AC\% and computational efficiency. Other values, such as 0.9 and 1.2, can also be selected as a slight difference is found in the solution quality compared to that of 0.6. When $\gamma$ increases to 0.4, the problem becomes more constrained, and the algorithm with $\theta_{max} = 0.6$ cannot solve all the instances within ten runs. In this case, the algorithm with setting $\theta_{max} = 0.9$ still outperforms the algorithm with other $\theta_{max}$ settings in terms of solution quality. The problem is highly constrained when $\gamma = 0.7$, and some instances may not have feasible solutions among ten runs. From the results, $\theta_{max} = 1.8$ has the highest proportion of feasible solutions compared to the algorithm with other $\theta_{max}$ values. The DA algorithm with setting $\theta_{max} = 0.9$ has a number of feasible solutions slightly less than that of $\theta_{max} = 1.8$. From the overall performance, we conclude that $\theta_{max}= 0.9$ can provide us with good solution quality and acceptable computational time in all the cases. We set $\theta_{max}=0.9$ in all the further experiments. For values of $\Theta_{red}$ and $n_{imp}$, we keep the initial settings, i.e, $\Theta_{red}=300$ and $n_{imp}=50$. 

\subsubsection{Contribution of local search operators}

As the algorithm largely relies on local search operators, their usefulness is verified. In this part, we analyze the contribution of local search operators to improve the solution quality. The effectiveness of each local search operator is presented, and the results of six different algorithm configurations are shown in Table \ref{operators}. In each of these configurations, one operator is excluded from the algorithm, and we run each algorithm configuration ten times, with each run iterating the respective algorithm 10000 times. We calculate the average solution gap of $BC\%$, $AC\%$, $Q1\%$, and $Q3\% $. Results for different algorithm configurations setting the previously selected parameter values ($\theta_{max}=0.9$) are summarized in Table \ref{operators}. For the scenario $\gamma = 0.7$, we report CPU times and FeasRatio.
 
\begin{table}[h!]
\renewcommand\arraystretch{0.8}
\centering
\small
\begin{threeparttable}
\caption{Experimental results when removing a single operator: \textit{Ex-pickup}, \textit{Ex-dropoff}, \textit{Ex-2-neighbor}, 
\textit{Relocate}, \textit{Exchange}, and \textit{2-opt}}
    \label{operators}
    \setlength\tabcolsep{7.5pt}
    \begin{tabular}{c c c c c c c c}
    \toprule
    Removing &None & \textit{Ex-pickup} &\textit{Ex-dropoff} &\textit{Ex-2-neighbor} &\textit{Relocate} &\textit{Exchange} &\textit{2-opt}   \\
    \midrule
    \textbf{$\gamma = 0.1$} \\
    \cline{1-1}
    $BC\%$  &0.10\%  &0.14\%  & 0.23\%  & 0.19\%  & 0.25\%  & 0.38\%  & 2.64\% \\
    $AC\%$  &0.52\%  &0.52\%  & 0.55\%  & 0.56\%  & 1.16\%  & 0.68\%  & 5.60\%\\
    $Q1\%$  &0.30\%  &0.40\%   & 0.40\%   & 0.44\%  & 0.79\%  & 0.51\%  & 3.76\%\\
    $Q3\%$  &0.73\%  &0.74\%  & 0.90\%   & 0.79\%  & 1.64\%  & 1.00\%   & 6.19\%\\
    FeasRatio &140/140 &140/140 &140/140 &140/140 &139/140 &140/140 &140/140\\
    CPU (s) &77.43  &74.88 & 71.41 & 88.97 & 57.53 & 79.51 & 68.92\\
    \hline
    \textbf{$\gamma = 0.4$} \\
    \cline{1-1}
    $BC\%$ &0.27\%  &0.27\%   & 0.27\%   & 0.38\%   & 0.38\%  & 0.27\%   & 2.56\%\\
    $AC\%$ &0.68\%  &0.73\%   & 0.74\%   & 0.78\%   & 1.15\%  & 0.84\%   & 4.92\%  \\
    $Q1\%$ &0.49\%  &0.51\%   & 0.49\%   & 0.64\%   & 0.86\%  & 0.63\%   & 3.66\% \\
    $Q3\%$ &0.84\%  &0.93\%   & 1.22\%   & 1.06\%   & 1.63\%  & 1.10\%    & 6.03\%\\
    FeasRatio &140/140  &140/140  & 140/140  & 139/140  & 136/140 & 140/140  & 140/140\\
    CPU (s) &116.97  &109.24 & 106.25 & 134.29 & 81.92 & 115.52 & 105.08\\
     \hline
    \textbf{$\gamma = 0.7$} \\
    \cline{1-1}
    FeasRatio &106/140  &96/140   & 106/140  & 90/140  & 86/140   & 97/140   & 74/140\\
    CPU (s) &201.68  &191.54 & 185.69 & 237.5 & 137.17 & 210.65 & 182.31\\
    \bottomrule
    \end{tabular}
\end{threeparttable}
\end{table}

We can find that each operator performs very well in improving the solution quality, especially the \textit{2-opt} operator. Additionally, the \textit{relocate} and \textit{2-opt} operator contributes to provide more feasible solutions in the case of $\gamma = 0.4, 0.7$. Therefore, it is necessary to include these operators in local search. As for \textit{add-request}, it is essential for inserting requests that are not served in the current solution. From the above analysis, the usefulness of each local search operator is proved.

 \subsubsection{Sensitivity analysis on number of iterations}

Then, we conduct the sensitivity analysis for the number of iterations $N_{iter}$. To identify a good $N_{iter}$, we conduct experiments with all the energy-level restrictions on type-a instances. We test ten values of $N_{iter}$, and report $BC\%$, $AC\%$, $Q1\%$, $Q3\%$ over ten runs. For the scenario of $\gamma = 0.7$, as different settings of $N_{iter}$ result in a different number of feasible solutions, we compare FeasRatio. The results are shown in Table \ref{dispersion alg}. 

\begin{table}[ht]
 \renewcommand\arraystretch{0.8}
\centering
 \small
\begin{threeparttable}
    \caption{Statistical comparison of DA performance under different iteration times for all $\gamma$ values on type-a instances}
    \label{dispersion alg}
    \setlength\tabcolsep{6pt}
    \begin{tabular}{c c c c c c c c c c c}
    \toprule
    $N_{iter}$ & 1000& 2000 & 3000& 4000& 5000& 6000 &7000 &8000 &9000 &10000 \\
    \hline
    \multicolumn{11}{c}{\textbf{Low energy restriction $\gamma = 0.1$}}\\
    \hline
    $BC\%$  &0.60\%   & 0.44\%  & 0.35\%  & 0.31\%  & 0.20\%   & 0.17\%  & 0.14\%  & 0.13\%  & 0.11\%  & 0.10\% \\
    $AC\%$   & NC   & NC   & 0.95\%  & 0.82\%  & 0.73\%  & 0.68\%  & 0.63\%  & 0.59\%  & 0.56\%  & 0.53\%\\
    $Q1\%$   & 1.42\%  & 0.82\%  & 0.61\%  & 0.52\%  & 0.45\%  & 0.42\%  & 0.36\%  & 0.34\%  & 0.31\%  & 0.30\% \\
    $Q3\%$   & 2.35\%  & 1.69\%  & 1.12\%  & 1.02\%  & 0.96\%  & 0.88\%  & 0.83\%  & 0.79\%  & 0.76\%  & 0.73\%\\
    FeasRatio &138/140 & 139/140 & 140/140 & 140/140 & 140/140 & 140/140 & 140/140 & 140/140 & 140/140 & 140/140\\
    CPU (s) & 10.15 & 17.5  & 24.86 & 32.43 & 39.92 & 47.45 & 54.92 & 62.38 & 69.79 & 77.43\\
    \hline
    \multicolumn{11}{c}{\textbf{Medium energy restriction $\gamma = 0.4$}}\\
    \hline
    $BC\%$   &1.07\%  & 0.72\%  & 0.57\%  & 0.48\%  & 0.40\%   & 0.37\%  & 0.34\%  & 0.31\%  & 0.30\%   & 0.27\%   \\
    $AC\%$   & NC   & NC   & NC   & 1.17\%  & 1.03\%  & 0.90\%   & 0.84\%  & 0.78\%  & 0.73\%  & 0.68\%\\
    $Q1\%$   & 1.69\%  & 1.18\%  & 0.96\%  & 0.82\%  & 0.72\%  & 0.66\%  & 0.63\%  & 0.57\%  & 0.54\%  & 0.49\%\\
    $Q3\%$   & 2.98\%  & 2.09\%  & 1.61\%  & 1.39\%  & 1.22\%  & 1.14\%  & 1.08\%  & 0.99\%  & 0.87\%  & 0.84\%\\
   FeasRatio &138/140 & 139/140 & 139/140 & 140/140 & 140/140 & 140/140 & 140/140 & 140/140 & 140/140 & 140/140\\
    CPU (s) &14.45 & 25.82 & 37.21 & 48.62 & 59.94 & 71.29 & 82.74 & 94.12 & 105.7 & 116.97 \\
    \hline
    \multicolumn{11}{c}{\textbf{High energy restriction $\gamma = 0.7$}}\\
    \hline
    FeasRatio &79/140 & 88/140 & 94/140 & 95/140 & 96/140 & 97/140 & 100/140 & 102/140 & 103/140 & 106/140 \\ 
    CPU (s) & 21.94 & 41.83 & 61.7 & 81.88 & 101.73 & 121.63 & 141.56 & 161.64 & 181.63 & 201.68 \\
    \bottomrule
    \end{tabular}
   \end{threeparttable}
\end{table}

From Table \ref{dispersion alg}, we observe that the values of $BC\%$, $AC\%$, $Q1\%$, $Q3\%$ are improved with more iterations. Among ten values of $N_{iter}$, 10000 iterations provide us with the best solution quality. We therefore set $N_{iter}$ to 10000 to conduct experiments. The performance of DA is also demonstrated as small results dispersion is found under all the values of $N_{iter}$. Moreover, we also notice that the computational time grows approximately linearly with the number of iterations, which is a computational advantage compared with the B\&C algorithm.

Note that choosing $N_{iter}=8000$ or $N_{iter}=9000$ slightly degrades the performances.
With such parameters, the computational time will be decreased.
Choosing $N_{iter}=10000$ is more robust, especially keeping in mind the evaluation of larger type-r instances.

\subsection{DA algorithm performance on the E-ADARP instances}
\label{DA performance}

In this section, we present the performance of our DA algorithm after tuning parameters from the previous section. Table \ref{cordeau instances results}, Table \ref{uber instances results}, and Table \ref{ropke instances results} present our DA algorithm results on type-a, -u, and -r instances under $\gamma = 0.1,0.4,0.7$, respectively. In each table, we report the values of $BC$, $AC$, $Q1$, $Q3$, and their corresponding gaps with $BC'$ (presented in the column named ``$BC'$"). If we obtain better solutions than the best-reported results of \cite{bongiovanni2019electric}, we mark them in bold with an asterisk. We mark our solutions in bold if they are equal to those reported in \cite{bongiovanni2019electric}. 

It should be noted that we find strictly better integer solutions than the reported optimal results of \cite{bongiovanni2019electric} in case of $\gamma = 0.4, 0.7$. The reason is that in the model of \cite{bongiovanni2019electric}, the employed ``big M" values were not correctly computed. We refer to supplementary material for a more in-depth analysis and how the ``big M" values should be set correctly. To distinguish these incorrect results, we mark them in italics in the column of ``$BC'$" and mark our obtained solutions in bold with double stars. The corresponding $BC\%$ values are therefore negative.

\subsubsection{Type-a instances results under different energy restrictions}
We first conduct experiments on type-a instances considering different scenarios $\gamma = 0.1, 0.4, 0.7$. A higher $\gamma$ value means a higher minimum battery level that vehicles must keep when returning to the destination depot. Recalling that each recharging station can only be visited at most once. The E-ADARP model is more constrained with an increasing $\gamma$. In Table \ref{cordeau instances results}, we compare our algorithm results to the best reported results in \cite{bongiovanni2019electric}.

We obtain equal/improved solutions for 36 out of 42 instances. Among them, 13 are the new best solutions. For some instances, we obtain better solutions than the reported optimal solutions in \cite{bongiovanni2019electric}, leading to savings on solution cost with up to 1.40\%. These instances are: a2-24-0.4, a3-30-0.4, a3-36-0.4, a2-24-0.7, a3-24-0.7, and a4-24-0.7. 

In all the scenarios, the proposed DA algorithm has quite small gaps to the best-reported results in \cite{bongiovanni2019electric}. In case of $\gamma = 0.1,0.4$, the average $BC\%$ is 0.05\% and 0.13\% (the worst $BC\%$ is 0.40\% and 1.26\%), and other values $AC\%$, $Q1\%$, $Q3\%$ are under 2.09\%, 2.51\%, and 2.79\%, respectively. 
When $\gamma = 0.7$, we consistently provide new solutions for a2-20, a4-32, and a5-40, while B\&C cannot solve these instances optimally or feasibly within two hours. Particularly, the generated new solutions on instance a4-32 and a5-40 have a much lower solution cost compared to the former reported best solutions in \cite{bongiovanni2019electric}, with an average gap of -7.49
\% and -5.22\%, respectively. 

In terms of computational efficiency, the CPU time for the proposed DA algorithm grows approximately in a linear way with sizes of instances. The average CPU time for all instances is 96.71s, and the proposed DA algorithm can efficiently solve large-scale instances within maters of minutes. 

\begin{table}[!hp]
 \renewcommand\arraystretch{0.8}
 \caption{Results of the proposed DA algorithm on type-a instances under $\gamma = 0.1, 0.4, 0.7$}
 \label{cordeau instances results}
    \begin{center}
    \footnotesize
    \begin{tabular}{c c c c c c c c c c c c}
    \hline
        \textbf{$\gamma = 0.1$}&\multicolumn{9}{c}{Proposed DA algorithm, 10000 iterations, 50 runs} & \multicolumn{2}{c}{Bongiovanni et al.,$^a$}\\
        \hline
        Instance & $BC$ & $BC\%$ & $Q1$ & $Q1\%$ & $AC$ &$AC\%$ &$Q3$  & $Q3\%$ &CPU(s) & $BC'$ &CPU$'$(s)\\
        \hline
        a2-16 &\textbf{237.38} & 0 & 237.38 & 0 & 237.38 & 0 & 237.38 & 0 & 39.3 & \textbf{237.38$^*$}  & 1.2  \\
        a2-20 &\textbf{279.08} & 0 & 279.08 & 0 & 279.08 & 0 & 279.08 & 0 & 73.8  &\textbf{279.08$^*$}  &4.2  \\
        a2-24 &\textbf{346.21} & 0  & 346.21 & 0  & 346.21 & 0 & 346.21 & 0  & 160.6 &\textbf{346.21$^*$} &9.0 \\
        a3-18 &\textbf{236.82} & 0 & 236.82 & 0 & 236.82 & 0 & 236.82 & 0 & 25.2  &\textbf{236.82$^*$}  &4.8  \\
        a3-24 &\textbf{274.80} & 0 & 274.80 & 0 & 274.80 & 0 & 274.80 & 0 & 58.3  &\textbf{274.80$^*$}  &13.8\\
        a3-30 &\textbf{413.27} & 0 & 413.27 & 0 & 413.27 & 0 & 413.27 & 0 & 54.3  &\textbf{413.27$^*$} &102.0\\
        a3-36 &\textbf{481.17} & 0  & 481.17 & 0  & 481.17 & 0  & 481.17 & 0  & 152.5 &\textbf{481.17$^*$} &106.8  \\
        a4-16 &\textbf{222.49} & 0 & 222.49 & 0 & 222.49 & 0 & 222.49 & 0 & 19.5 &\textbf{222.49$^*$}  & 3.6 \\
        a4-24 &\textbf{310.84} & 0 & 310.84 & 0 & 310.84 & 0 & 312.44 & 0.51\%  & 29.6  &\textbf{310.84$^*$} &31.2\\
        a4-32 &\textbf{393.96} & 0 & 393.95 & 0 & 395.12 & 0.29\%  & 397.58 & 0.92\%  & 52.0 &\textbf{393.96$^*$} &612.0 \\
        a4-40 &\textbf{453.84} & 0 & 458.22 & 0.97\%  & 459.42 & 1.23\%  & 460.56 & 1.48\%  & 92.0  &\textbf{453.84$^*$} &517.2 \\
        a4-48 &555.93 & 0.25\%  & 560.19 & 1.02\%  & 561.26 & 1.21\%  & 562.87 & 1.50\%  & 141.8 &\textbf{554.54} & 7200.0  \\
        a5-40 &414.80 & 0.07\%  & 418.48 & 0.96\%  & 420.35 & 1.41\%  & 422.56 & 1.94\%  & 64.9  &\textbf{414.51$^*$}  &1141.8 \\
        a5-50 &561.41 & 0.40\%  & 567.82 & 1.55\%  & 570.58 & 2.04\%  & 573.51 & 2.56\%  & 137.3   &\textbf{559.17} &7200.0\\
        \hline
        Summary & & 0.05\%  &   & 0.32\%  &  & 0.44\%  &  & 0.64\%  & 78.6 & &1210.5 \\
        \hline
        \textbf{$\gamma = 0.4$} & $BC$ & $BC\%$ & $Q1$ & $Q1\%$ & $AC$ &$AC\%$ &$Q3$  & $Q3\%$ &CPU(s) & $BC'$ &CPU$'$(s)\\
        \hline
        a2-16 &\textbf{237.38} & 0 & 237.38 & 0 & 237.38 & 0 & 237.38 & 0 & 52.9 &\textbf{237.38$^*$}  &1.8\\
        a2-20 &\textbf{280.70} & 0  & 280.70 & 0 & 280.70 & 0 & 280.70 & 0  & 140.7 &\textbf{280.70$^*$}  &49.8 \\
        a2-24 &\textbf{347.04$^{**}$} & -0.29\% & 347.04 & -0.29\% & 347.04 & -0.29\% & 347.04 & -0.29\% & 231.0 &\textit{348.04$^*$} &25.2 \\
        a3-18 &\textbf{236.82} & 0 & 236.82 & 0 & 236.82 & 0 & 236.82 & 0 & 26.3  &\textbf{236.82$^*$}  &4.2 \\
        a3-24 &\textbf{274.80} & 0  & 274.80 & 0  & 274.80 & 0  & 276.11 & 0.48\%  & 67.9 &\textbf{274.80$^*$}  &16.8\\
        a3-30 &\textbf{413.34$^{**}$} & -0.01\% & 413.34 & -0.01\% & 413.34 & -0.01\% & 413.34 & -0.01\% & 88.7 &\textit{413.37$^*$}  &99.0 \\
        a3-36 &\textbf{483.06$^{**}$} & -0.22\% & 483.83 & -0.06\% & 483.86 & -0.06\% & 485.43 & 0.27\%  & 157.8 &\textit{484.14$^*$} &306.6 \\
        a4-16 &\textbf{222.49} & 0 & 222.49 & 0 & 222.49 & 0 & 222.49 & 0 & 19.4 &\textbf{222.49$^*$} &5.4 \\
        a4-24 &\textbf{311.03} & 0  & 311.28 & 0.08\%  & 311.65 & 0.20\%  & 313.21 & 0.70\%  & 32.0 &\textbf{311.03$^*$} &39.6 \\
        a4-32 &\textbf{394.26} & 0 & 395.05 & 0.20\%  & 397.21 & 0.75\%  & 400.32 & 1.54\%  & 63.0   &\textbf{394.26$^*$}  &681.6 \\
        a4-40 &\textbf{453.84} &0 & 457.20 & 0.74\%  & 459.46 & 1.24\%  & 461.06 & 1.59\%  & 116.7  &\textbf{453.84$^*$}  &417.6 \\
        a4-48 &558.11 & 0.63\%  & 561.40 & 1.23\%  & 563.47 & 1.60\%  & 565.35 & 1.94\%  & 177.5 &\textbf{554.60}  &7200.0 \\
        a5-40 &416.25 & 0.42\%  & 418.97 & 1.08\%  & 420.32 & 1.40\%  & 422.75 & 1.99\%  & 72.6  &\textbf{414.51$^*$}  &1221.0 \\
        a5-50 &567.54 & 1.26\%  & 572.23 & 2.09\%  & 574.56 & 2.51\%  & 576.11 & 2.79\%  & 162.8  &\textbf{560.50}  &7200.0 \\
        \hline
        Summary & & 0.13\%  &  & 0.36\%  &   & 0.52\%  & & 0.79\%  & 100.7 & &1233.5 \\
        \hline
        \textbf{$\gamma = 0.7$} & $BC$ & $BC\%$ & $Q1$ & $Q1\%$ & $AC$ &$AC\%$ &$Q3$  & $Q3\%$ &CPU(s) & $BC'$ &CPU$'$(s)\\
        \hline
        a2-16 &\textbf{240.66} & 0  & 240.66 & 0  & 240.66 & 0 & 240.66 & 0  & 95.8 &\textbf{240.66$^*$}  &5.4 \\
        a2-20 &\textbf{293.27$^*$} & --   & 293.27 & --   & 294.11 & --   & NA & NA   & 172.8 &NA  &7200.0  \\
        a2-24 &\textbf{353.18$^{**}$} & -1.40\% & 366.49 & 2.31\%  & NA    & NA   & NA    & NA  & 206.6 &\textit{358.21$^*$}  &961.2  \\
        a3-18 &\textbf{240.58} & 0 & 240.58 & 0 & 240.58 & 0 & 240.58 & 0 & 58.3  &\textbf{240.58$^*$}  &48.0 \\
        a3-24 &\textbf{275.97$^{**}$} & -0.63\% & 275.97 & -0.63\% & 277.43 & -0.10\% & 279.13 & 0.51\%  & 123.7  &\textit{277.72$^*$}  &152.4 \\
        a3-30 &\textbf{424.93$^*$} & --    & 432.29 & --   & 436.20 & --   & NA    & NA   & 77.7   &NA  &7200.0 \\
        a3-36 &\textbf{494.04} & 0 & 497.11 & 0.62\%  & 502.27 & 1.67\%  & 505.95 & 2.41\%  & 125.4  &\textbf{494.04}  &7200.0 \\
        a4-16 &\textbf{223.13} & 0  & 223.13 & 0 & 223.13 & 0  & 223.13 & 0  & 31.3 &\textbf{223.13$^*$} &67.2 \\
        a4-24 &\textbf{316.65$^{**}$} & -0.49\% & 318.21 & 0 & 318.31 & 0.03\%  & 320.87 & 0.84\%  & 53.7   &\textit{318.21$^*$}  &1834.8 \\
        a4-32 &\textbf{397.87$^*$} & -7.49\% & 401.58 & -6.63\% & 405.85 & -5.63\% & 408.69 & -4.97\% & 71.4  &430.07  &7200.0 \\
        a4-40 &\textbf{479.02$^*$} & --   & NA    & NA   & NA    & NA   & NA    & NA   & 114.7 &NA  &7200.0 \\
        a4-48 &\textbf{582.22$^*$}   & --   &610.75   & --   & NA    & NA   & NA    & NA   & 164.4 &NA  &7200.0 \\
        a5-40 &\textbf{424.26$^*$} & -5.22\% & 433.12 & -3.24\% & 436.94 & -2.39\% & 441.15 & -1.45\% & 97.5 &447.63  &7200.0\\
        a5-50 &\textbf{603.24$^*$}    & --   & NA    & NA   & NA    & NA   & NA    & NA   & 158.4 &NA &7200.0 \\
        \hline
        Summary & &-- & &-- & &-- & &-- &110.8 & &4333.4 \\
        \hline
    \end{tabular}
    \begin{tablenotes}
    \footnotesize
     \item a: Due to incorrect big M values, some of the reported optimal results of \cite{bongiovanni2019electric} are higher than our obtained solution values. Those results are highlighted in italics and our obtained results are marked in bold with double stars;
   \end{tablenotes}
    \end{center}
\end{table}

\begin{table}[!hp]
\renewcommand\arraystretch{0.8}
\caption{Results of the proposed DA algorithm on type-u instances under $\gamma = 0.1, 0.4, 0.7$}
\label{uber instances results}
    \begin{center}
    \footnotesize
    \begin{tabular}{c c c c c c c c c c c c}
    \hline
        \textbf{$\gamma = 0.1$}&\multicolumn{9}{c}{Proposed DA algorithm, 10000 iterations, 50 runs} & \multicolumn{2}{c}{Bongiovanni et al.,$^a$}\\
        \hline
        Instance & $BC$ & $BC\%$ & $Q1$ & $Q1\%$ & $AC$ &$AC\%$ &$Q3$  & $Q3\%$ &CPU(s) & $BC'$ &CPU$'$(s)\\
        \hline
        u2-16 &\textbf{57.61}  & 0  & 57.61  & 0  & 57.61  & 0  & 57.61  & 0  & 120.1 & \textbf{57.61$^*$}  & 21.0  \\
        u2-20 &\textbf{55.59}  & 0 & 55.59  & 0 & 56.34  & 1.34\%  & 56.34  & 1.34\%  & 401.8  &\textbf{55.59$^*$}  &9.6  \\
        u2-24 &\textbf{90.73$^{**}$}  & -0.60\%  & 90.84  & -0.47\% & 90.84  & -0.47\% & 90.98  & -0.32\% & 599.7 & \textit{91.27$^*$} &432.0 \\
        u3-18 &\textbf{50.74}  & 0  & 50.74  & 0  & 50.74  & 0  & 50.93  & 0.37\%  & 108.3 &\textbf{50.74$^*$}  &10.8  \\
        u3-24 &\textbf{67.56}  & 0 & 67.87  & 0.46\%  & 68.16  & 0.89\%  & 68.16  & 0.89\%  & 111.5 &\textbf{67.56$^*$}  &130.2\\
        u3-30 &\textbf{76.75}  & 0 & 77.21  & 0.60\%  & 77.80  & 1.37\%  & 78.65  & 2.47\%  & 174.1 &\textbf{76.75$^*$} &438.0\\
        u3-36 &104.27 & 0.22\%  & 104.87 & 0.79\%  & 105.42 & 1.33\%  & 106.36 & 2.23\%  & 420.7 &\textbf{104.04$^*$} &1084.8  \\
        u4-16 &\textbf{53.58}  & 0  & 53.58  &0 & 53.58  & 0  & 53.58  & 0 & 51.4 &\textbf{53.58$^*$}  & 48.0 \\
        u4-24 &90.13  & 0.34\%  & 90.72  & 1.00\%  & 90.85  & 1.14\%  & 90.95  & 1.25\%  & 55.3 &\textbf{89.83$^*$} &13.2\\
        u4-32 &\textbf{99.29}  & 0  & 99.29  & 0  & 99.42  & 0.13\%  & 99.67  & 0.38\%  & 119.1  &\textbf{99.29$^*$} &1158.6 \\
        u4-40 &\textbf{133.11} & 0 & 134.46 & 1.02\%  & 135.18 & 1.55\%  & 136.08 & 2.23\%  & 154.0 &\textbf{133.11$^*$} &185.4 \\
        u4-48 &\textbf{147.75$^*$} & -0.37\% & 148.87 & 0.39\%  & 149.69 & 0.93\%  & 150.42 & 1.43\%  & 841.0 &148.30 & 7200.0  \\
        u5-40 &\textbf{121.86} & 0 & 123.11 & 1.03\%  & 123.38 & 1.25\%  & 124.47 & 2.14\%  & 113.8  &\textbf{121.86}  &1141.8 \\
        u5-50 &144.22 & 0.78\%  & 145.04 & 1.36\%  & 145.63 & 1.77\%  & 146.30 & 2.24\%  & 245.5 &\textbf{143.10} &7200.0\\
        \hline
        Summary & & 0.03\%  &   & 0.44\%  &  & 0.80\%  & & 1.19\%  & 251.2 & &1795.1 \\
        \hline
        \textbf{$\gamma = 0.4$} & $BC$ & $BC\%$ & $Q1$ & $Q1\%$ & $AC$ &$AC\%$ &$Q3$  & $Q3\%$ &CPU(s) & $BC'$ &CPU$'$(s)\\
        \hline
        u2-16 &\textbf{57.65}  & 0 & 57.65  & 0 & 57.65  & 0 & 57.65  & 0 & 156.6 &\textbf{57.65$^*$}  &25.8\\
        u2-20 &\textbf{56.34}  & 0 & 56.34  & 0 & 56.34  & 0 & 56.34  & 0 & 606.6 &\textbf{56.34$^*$}  &12.0 \\
        u2-24 &\textbf{91.24$^{**}$}  & -0.43\% & 91.27  & -0.39\% & 91.72  & 0.10\%  & 92.06  & 0.47\%  & 817.8 &\textit{91.63$^*$} &757.2 \\
        u3-18 &\textbf{50.74}  & 0   & 50.74  & 0  & 50.74  & 0  & 50.99  & 0.50\%  & 125.0   &\textbf{50.74$^*$}  &13.8 \\
        u3-24 &\textbf{67.56}  & 0  & 67.87  & 0.46\%  & 68.16  & 0.89\%  & 68.16  & 0.89\%  & 141.0  &\textbf{67.56$^*$}  &220.8\\
        u3-30 &\textbf{76.75}  & 0 & 77.12  & 0.48\%  & 77.93  & 1.54\%  & 78.65  & 2.48\%  & 285.8   &\textbf{76.75$^*$}  &336.6 \\
        u3-36 &104.49 & 0.41\%  & 105.65 & 1.53\%  & 106.37 & 2.22\%  & 107.19 & 3.01\%  & 898.9  &\textbf{104.06$^*$} &2010.0 \\
        u4-16 &\textbf{53.58}  & 0  & 53.58  & 0  & 53.58  & 0  & 53.58  & 0  & 60.5 &\textbf{53.58$^*$} &44.4 \\
        u4-24 &90.72  & 1.00\%  & 90.72  & 1.00\%  & 91.00  & 1.30\%  & 91.12  & 1.44\%  & 65.6   &\textbf{89.83$^*$} &28.2 \\
        u4-32 &\textbf{99.29}  & 0 & 99.29  & 0  & 99.42  & 0.13\%  & 99.90  & 0.61\%  & 156.3   &\textbf{99.29$^*$}  &2667.6 \\
        u4-40 &\textbf{133.78$^{**}$} & -0.10\% & 135.43 & 1.14\%  & 135.83 & 1.44\%  & 136.56 & 1.98\%  & 303.1 &\textit{133.91$^*$}  &2653.2 \\
        u4-48 &\textbf{148.48$^*$} & --   & 149.86 & --   & 150.81 & --   & 151.77 & --   & 1390.7 &NA  &7200.0 \\
        u5-40 &\textbf{121.96$^*$} & -0.22\% & 123.08 & 0.69\%  & 123.63 & 1.15\%  & 124.42 & 1.79\%  & 160.8  &122.23  &7200.0 \\
        u5-50 &143.68 & 0.38\%  & 145.66 & 1.76\%  & 146.60 & 2.42\%  & 147.15 & 2.80\%  & 391.5 &\textbf{143.14}  &7200.0 \\
        \hline
        Summary &  & --  &   & --  &   & --  &  & -- & 397.2  & &2169.3 \\
        \hline
        \textbf{$\gamma = 0.7$} & $BC$ & $BC\%$ & $Q1$ & $Q1\%$ & $AC$ &$AC\%$ &$Q3$  & $Q3\%$ &CPU(s) & $BC'$ &CPU$'$(s)\\
        \hline
        u2-16 &\textbf{59.19}  & 0  & 59.26  & 0.11  & 60.01  & 1.38  & 60.19  & 1.69  & 419.6  &\textbf{59.19$^*$}  &338.4 \\
        u2-20 &\textbf{56.86}  & 0  & 58.39  & 2.69  & 58.39  & 2.69  & 58.88  & 3.55  & 1527.6 &\textbf{56.86$^*$}  &72.0  \\
        u2-24 &\textbf{92.84$^*$}  & --  & 94.33    & --  & 99.38    & --   & NA   & NA  & 502.5 &NA  &7200.0  \\
        u3-18 &\textbf{50.99}  & 0  & 50.99  & 0 & 50.99  & 0  & 50.99 & 0 & 206.9   &\textbf{50.99$^*$}  &24.0 \\
        u3-24 &\textbf{68.39}  & 0 & 68.39  & 0 & 68.44  & 0.08\%  & 68.73  & 0.49\%  & 375.8  &\textbf{68.39$^*$}  &400.2 \\
        u3-30 &\textbf{77.94$^{**}$}  & -0.26\% & 78.72  & 0.74\%  & 79.37  & 1.57\%  & 79.56  & 1.81\%  & 1094.8 &\textit{78.14$^*$}  &3401.4 \\
        u3-36 &106.00 & 0.20\%  & 106.41 & 0.59\%  & 107.57 & 1.68\%  & 107.92 & 2.01\%  & 1606.4  &\textbf{105.79}  &7200.0 \\
        u4-16 &\textbf{53.87}  & 0 & 53.87  & 0 & 53.87  & 0 & 53.87  & 0 & 96.9 &\textbf{53.87$^*$}  &88.8 \\
        u4-24 &90.07  & 0.12\%  & 90.97  & 1.12\%  & 90.97  & 1.12\%  & 90.97  & 1.12\%  & 254.5   &\textbf{89.96$^*$}  &22.8 \\
        u4-32 &\textbf{99.50}  & 0 & 100.01 & 0.51\%  & 101.09 & 1.60\%  & 101.75 & 2.26\%  & 325.3 &\textbf{99.50$^*$}  &2827.2 \\
        u4-40 &\textbf{136.08$^*$} & --   & 137.65 & --   & 138.98 & --   & NA    & NA   & 708.0  &NA  &7200.0 \\
        u4-48  &\textbf{152.58$^*$} & --   & 157.85 & --   & 162.62 & --   & NA    & NA   & 1958.8  &NA  &7200.0 \\
        u5-40 &\textbf{123.52$^*$} & --  & 125.30 & --   & 126.10 & --   & 127.08 & --   & 359.6 &NA &7200.0\\
        u5-50 &\textbf{143.51$^*$} & -0.59\% & 148.16 & 2.64\%  & 149.52 & 3.58\%  & 152.36 & 5.54\%  & 922.2  &144.36 &7200.0 \\
        \hline
        Summary & &-- & &-- & &-- & &-- &780.1 & &3598.2 \\
        \hline
    \end{tabular}
    \begin{tablenotes}
    \footnotesize
     \item a: Due to incorrect big M values, some of the reported optimal results of \cite{bongiovanni2019electric} are higher than our obtained solution values. Those results are highlighted in italics and our obtained results are marked in bold with double stars;
   \end{tablenotes}
    \end{center}
\end{table}

\subsubsection{Type-u instances results under different energy restrictions}

On type-u instances, we conduct experiments under different energy-restriction levels $\gamma = 0.1, 0.4, 0.7$. The results are shown in Table \ref{uber instances results}.

The proposed DA algorithm finds equal solutions for 22 out of 42 instances and finds new best solutions for 12 previously solved and unsolved instances. Particularly, on instance u2-24-0.1, u2-24-0.4, u4-40-0.4, and u3-30-0.7, we find strictly better solutions than the reported optimal solutions in \cite{bongiovanni2019electric}, which contributes to savings on solution costs with up to 0.43\%. In each scenario, our best solutions have quite small gaps to the $BC'$ reported in \cite{bongiovanni2019electric} (the worst- and best-case $BC\%$is 1.00\% and -0.60\%, respectively). We further demonstrate our algorithm consistency via other statistical values ($Q1\%$, $AC\%$, $Q3\%$), as our algorithm continuously finds high-quality solutions with the increasing size of instances. 
In terms of computational efficiency, solving the problem exactly seems more computationally effective on small-sized instances. The reason is that we fix $N_{iter}$ to 10000 for all instances, whereas the small-sized ones can be solved to their best values (i.e., optimal objective values reported in \cite{bongiovanni2019electric}) in much fewer iterations. 
On medium-to-large-sized instances, using an efficient heuristic (e.g., the proposed DA algorithm) is a more computational appealing option.


\subsubsection{Type-r instances results under different energy restrictions}
\label{large-scale type-r}
We present our algorithm results on type-r instances in Table \ref{ropke instances results}. These results are the first solutions found for these new instances and can serve as benchmark results for future studies.

In scenarios $\gamma = 0.1$ and $\gamma = 0.4$, we find feasible solutions for 19 out of 20 instances, with an average CPU time of 269.71s and 373.89s, respectively. 
When increasing from $\gamma=0.1$ to $\gamma=0.4$, the statistical dispersion also increases, but the dispersion remains quite acceptable. For instance r7-84, most of the runs with $\gamma=0.4$ do not find a feasible solution. For instance r8-96, our DA algorithm cannot find a feasible solution among 50 runs with $\gamma=0.4$. These instances seem challenging for future works.

When $\gamma = 0.7$, we found no feasible solution for all the type-r instances, despite 50 runs and 10000 iterations. One reason is that many of these instances are too constrained to be feasible for $\gamma=0.7$ with the limitation of visiting recharging stations. However, it opens a perspective to prove it using exact methods with lower bounds. 

\begin{table}[h!]
\renewcommand\arraystretch{0.6}
    \centering
     \small
    \begin{threeparttable}
    \caption{Results of the proposed DA algorithm with 10000 iterations 50 runs on type-r instances under $\gamma = 0.1, 0.4$}
    \label{ropke instances results}
    \setlength\tabcolsep{7.5pt}
    \begin{tabular}{c c c c c c c c c}
        \toprule
        \textbf{$\gamma = 0.1$} & $BC$ &$Q1$ & $Q1\%$  &$AC$ &$AC\%$  &$Q3$ & $Q3\%$ &CPU(s)\\
        \midrule
        r5-60 &691.83  & 699.93  & 1.17\% & 706.20  & 2.08\% & 710.43  & 2.69\% & 178.44 \\
        r6-48 &506.72  & 509.67  & 0.58\% & 512.69  & 1.18\% & 515.39  & 1.71\% & 229.31 \\
        r6-60 &692.00  & 696.67  & 0.67\% & 700.15  & 1.18\% & 703.95  & 1.73\% & 127.03  \\
        r6-72 &777.44  & 788.12  & 1.37\% & 794.69  & 2.22\% & 801.87  & 3.14\% & 208.39 \\
        r7-56 &613.10  & 620.69  & 1.24\% & 624.51  & 1.86\% & 630.72  & 2.87\% & 88.20  \\
        r7-70 &760.90  & 772.45  & 1.52\% & 778.84  & 2.36\% & 786.02  & 3.30\% & 209.76  \\
        r7-84 &889.38  & 900.34  & 1.23\% & 904.88  & 1.74\% & 913.88  & 2.75\% & 322.66 \\
        r8-64 &641.99  & 647.87  & 0.92\% & 652.59  & 1.65\% & 657.49  & 2.41\% & 612.06\\
        r8-80 &803.52  & 820.96  & 2.17\% & 828.67  & 3.13\% & 834.19  & 3.82\% & 357.75 \\
        r8-96 &1053.11 & 1069.98 & 1.60\% & 1080.80 & 2.63\% & 1089.96 & 3.50\% & 363.46  \\
        \hline
        Summary &&  & 1.25\%   &   & 2.00\%  & &2.79\% &269.71\\
        \hline
        \textbf{$\gamma = 0.4$} & $BC$ &$Q1$ & $Q1\%$  &$AC$ &$AC\%$  &$Q3$ & $Q3\%$ &CPU(s)\\
        \hline
        r5-60 &697.97 & 710.30 & 1.77\% & 718.44 & 2.93\% & 727.27 & 4.20\% & 293.25 \\
        r6-48 &506.91 & 509.48 & 0.51\% & 514.46 & 1.49\% & 517.53 & 2.10\% & 257.59 \\
        r6-60 &694.78 & 702.67 & 1.14\% & 706.07 & 1.62\% & 710.80 & 2.31\% & 173.43 \\
        r6-72 &799.60 & 811.85 & 1.53\% & 821.17 & 2.70\% & 832.07 & 4.06\% & 349.98  \\
        r7-56 &613.66 & 620.58 & 1.13\% & 624.40 & 1.75\% & 627.51 & 2.26\% & 99.91 \\
        r7-70 &766.05 & 778.70 & 1.65\% & 784.54 & 2.41\% & 791.07 & 3.27\% & 273.52 \\
        r7-84 &932.12 & 964.04 & 3.43\% & NA   & NA  & NA   & NA  & 584.26 \\
        r8-64 &638.36 & 649.84 & 1.80\% & 652.30 & 2.18\% & 657.02 & 2.92\% & 641.63 \\
        r8-80 &811.19 & 823.70 & 1.54\% & 833.05 & 2.69\% & 841.76 & 3.77\% & 448.14\\
        r8-89 &NA & NA   & NA  & NA    & NA  & NA   & NA   & 617.17 \\
        \hline
        Summary & &  &NA &   &NA &  &NA  &373.89\\
        \bottomrule
    \end{tabular}
   \end{threeparttable}
   \vspace{-4mm}
\end{table}

\subsubsection{Conclusion of algorithm performance}

On both type-a and -u instances, we observe the limit of solving capabilities of the B\&C. Even with a time limit of two hours, it is difficult for B\&C to solve medium-to-large-sized E-ADARP instances, especially under a high energy restriction. Our DA algorithm can continuously provide high-quality solutions for highly constrained instances within a reasonable computational time.
We also show that our DA algorithm can tackle larger-sized instances with up to 8 vehicles and 96 requests.
Nineteen type-r instances for $\gamma=0.1$ and $\gamma=0.4$ are solved feasibly, and these results are the first solutions found for these new instances, which can serve as a benchmark for future studies. To conclude, the proposed DA algorithm remains highly effective and can provide optimal/near-optimal solutions even facing highly constrained instances. The proposed DA algorithm significantly outperforms the B\&C algorithm for medium-to-large-sized instances, and its consistency seems quite acceptable for such difficult instances.

\subsection{Sensitivity analysis of the maximum number of charging visits per station}
\label{multiple section}

 As discussed in Section \ref{sec::multiple?}, the hypothesis of visiting each recharging station at most once is not realistic. We adjust our DA algorithm as mentioned in Section \ref{adapt DA to allow multiple visits} to allow multiple visits to each recharging station. The adjusted DA algorithm is able to investigate the effect of increasing the value of $n_{as}$ on solution cost and feasibility. Recalling that we analyze four different cases: $n_{as} =1,2,3,\infty$. 
 
 For type-a instances, as in the scenario of $\gamma = 0.1$, we obtain optimal solutions for most of the instances, and other instances are solved without visiting recharging stations. Therefore, we focus on scenarios of $\gamma = 0.4, 0.7$ and analyze the effect of allowing multiple visits in these cases. For type-u and -r instances, we conduct experiments with adjusted DA algorithm with $n_{as} = 2,3,\infty$ under $\gamma \in \{0.1,0.4, 0.7\}$. The detailed results are presented in \ref{unlimited visits}. In Table \ref{mul cordeau} and \ref{mul ropke}, we compare DA algorithm results on each instance with setting $n_{as} = 1,2,3,\infty$ and we mark the best one(s) in bold. In Table \ref{mul uber}, we compare our algorithm results under each setting of $n_{as}$ with the reported results in \cite{bongiovanni2019electric}.  Improved solutions are marked in bold with an asterisk while equal solutions are marked in bold. In the column of $BC_{\infty}$, if the obtained solution is better than other solutions obtained under $n_{as} = 1,2,3$, we mark it in bold with double stars. 
 On each instance, the adjusted DA algorithm performs 50 runs with 10000 iterations per run. We report the maximum number of recharging visits experienced on a station (denoted as $N_{max}^s$) for the best-obtained solution under $n_{as} = \infty$ in the column named ``$N_{max}^s$". In addition, we also report the average number of visited recharging stations under setting $n_{as} = \infty$ in the column of ``$\overline{N_{avg}^s}$".

From these results, we observe that the previous difficulties for the DA algorithm to solve the E-ADARP instances are reduced considering multiple visits per station. The major findings are: (1) significant increases on $\overline{N_{avg}^s}$ are observed on all instances with increasing $\gamma$ value, especially on type-r instances, where the average value of $\overline{N_{avg}^s}$ is tripled when $\gamma$ changes to 0.7; 
(2) allowing multiple visits to each recharging station improves the solution quality as we found lower-cost solutions. Particularly, we obtain feasible solutions for all type-r instances under $\gamma = 0.7$ with $n_{as} = 3, \infty$, while no feasible solution is found with $n_{as} = 1$; (3) for type-a and -r instances, relaxing to $n_{as} = \infty$ seems to be more computationally attractive as it does not introduce additional computational time, compared to the results obtained by replicating recharging stations. For type-u instances, having a pre-calculated $n_{as}$ would be more computationally favorable; 
(4) on average, allowing at-most-two and -three visits per station slightly increases the computational time. Allowing at-most-three visits per station seems to strike a good balance between solution quality and computational time; (5) $n_{as} = 3$ seems to be a good upper bound for solving type-u instances allowing multiple recharging visits, while one needs to set $n_{as}$ to 4 and 7 for type-a and -r instances, respectively. A potential perspective from these results would be to investigate more realistic constraints, e.g., on the capacity of recharging stations, rather than limiting visits to recharging stations in the E-ADARP. Another direction for future studies is to design a heuristic to calculate the upper bound on the total number of recharging visits for a given route.

\section{Conclusions and Perspectives} \label{sec::conclusion}

This paper proposes an efficient DA algorithm to solve the E-ADARP, which aims to minimize a weighted-sum objective, including the total travel time and the total excess user ride time. To minimize the total excess user ride time, we propose a fragment-based representation of paths. A new method is developed upon this representation to calculate the minimum excess user ride time for a given route. Another challenge in solving the E-ADARP involves incorporating the partial recharging at recharging stations, which complicates the feasibility checking of a given route; to resolve this issue, we propose an exact route evaluation scheme of linear time complexity that can accurately handle the effect of allowing partial recharging and validate the feasibility of solutions. These two methods compose an exact and efficient optimization of excess user ride time for an E-ADARP route. To the best of our knowledge, this is the first time that total excess user ride time is optimized in an exact way for the E-ADARP. 

In computational experiments, we first prove the effectiveness and accuracy of our DA algorithm compared to the best-reported results of \cite{bongiovanni2019electric}. On 84 existing E-ADARP instances, our DA algorithm obtains equal solutions for 45 instances and provides better solutions on 25 instances. We also demonstrate that the proposed DA algorithm can consistently provide high-quality solutions in a short computational time. On the previously solved instances, the DA algorithm improves the solution quality by 0.16\% on average. On newly introduced large-scale E-ADARP instances, we provide new solutions for 19 instances. These results may serve as benchmark results for future studies. We then extend the E-ADARP model to allow unlimited visits to each recharging station. The previous difficulties for DA local search are lessened under this more realistic situation, and the results are less dispersed than the results of the at-most-one visit to each recharging station.
Our extension of the E-ADARP model thus offers a new perspective in proposing a more realistic constraint in the E-ADARP for recharging stations, e.g., considering capacity and scheduling constraints in recharging stations.



Our results offer other new perspectives for the E-ADARP in terms of algorithmic and modeling aspects. First, some instances remain unsolvable even after 50 independent runs of the DA algorithm. One reason may be that no feasible solution exists for these instances, which remain challenging for future studies using heuristic and exact methods. An interesting investigation would be examining the effects of more randomness in the algorithm, for example, considering a sequence of randomly ordered operators. 
 Second, the computational efficiency of our algorithm could be further improved by applying a more intelligent insertion strategy of recharging stations, adapting a parallel version of the DA algorithm, and designing stopping criteria to terminate the algorithm before completing all iterations, as in \cite{ropke2006adaptive}.
The E-ADARP could also be extended to consider user's inconvenience as a second objective, which helps understand the conflicting interests between service providers and users and provide a high-quality approximation of Pareto front for decision makers. The proposed excess user ride time optimization approach can also be adapted to solve the classical DARP in the context of multiple objectives, in which the total excess user ride time is minimized as a separate objective.
Another way that the E-ADARP model may be improved involves taking into account more real-life characteristics. For example, time-dependent travel times occur with traffic jams in peak hours. 
Relatedly, the static E-ADARP can be extended to dynamic E-ADARP, taking into account updates of requests during the day (e.g., new requests, cancellations, modifications).
Quick and efficient routing and scheduling heuristics for the dynamic E-ADARP as in \cite{bongiovanni2022machine} is crucial in such a context where metaheuristics also seem promising.

\ACKNOWLEDGMENT{The authors would like to thank Claudia Bongiovanni, Mor Kaspi, and Nikolas Geroliminis for kindly providing access to their implementation. This work is supported by the China Scholarship Council (CSC, grant number 201807000145) and by public funding within the scope of the French Program ``Investissements d’Avenir”.
}

\bibliographystyle{informs2014trsc} 
\bibliography{mybib.bib} 

\newpage
\begin{APPENDICES}
\section{Proof of Theorem 1 and Linear Programming Model} \label{proof of theorem}
We recall Theorem \ref{theorem1} and present the proof below.
\subsection{Proof of Theorem 1 and Examples}
\textbf{Theorem \ref{theorem1}:} If $\mathcal{R}$ is a feasible route and   $\mathcal{F}_1, \mathcal{F}_2, \cdots, \mathcal{F}_n$ are all the fragments on $\mathcal{R}$, then we have $EU_{min}(\mathcal{R}) = EU_{min}(\mathcal{F}_1) + EU_{min}(\mathcal{F}_2) + \cdots + EU_{min}(\mathcal{F}_n)$

\begin{proof}
In this proof, a schedule is called ``optimal" if it has minimal excess user ride time.

Assuming that $\mathcal{T} = [\cdots,T_v,\cdots]_{v\in \mathcal{R}}$ is an optimal schedule of route $\mathcal{R}$, $T_v$ is the service start time at node $v$, and the arrival time of node $v$ is: $arr_v = T_{v-1} + t_{v-1,v} + s_{v-1}$. To prove the theorem, it is enough to show that for each fragment $\mathcal{F}_i \subseteq \mathcal{R}$, the restricted schedule  $\mathcal{T}|_{\mathcal{F}_i} = [\cdots,T_v,\cdots]_{v\in \mathcal{F}_i}$ over $\mathcal{F}_i$ is also an optimal schedule for $\mathcal{F}_i$. To simplify the notation, we denote $\mathcal{T}|_{\mathcal{F}_i}$ as $\mathcal{T}_i$.  
Our proof consists of two different cases :
\begin{enumerate}
    \item $arr_v= T_v$ for all $v\in \mathcal{F}_i$ and node $v$ is not the start node of $\mathcal{F}_i$. In this case, the vehicle starts service at its arrival on each node in $\mathcal{F}_i$. Clearly,  $\mathcal{T}_i$ is also an optimal schedule over $\mathcal{F}_i$ as the waiting time on $\mathcal{F}_i$ is zero, proof is finished;
    \item $arr_v< T_v$ for some $v\in \mathcal{F}_i$ and node $v$ is not the start node of $\mathcal{F}_i$. In this case, waiting time is generated at some nodes.
    Let $v_1 \in \mathcal{F}_i$ be the first node such that $arr_{v_1}< T_{v_1}$ and $v_2 \in \mathcal{F}_i$ be the last node such that $arr_{v_2}< T_{v_2}$. Then we derive the following properties of $\mathcal{T}_i$:
    \begin{itemize}
        \item[(i)] $T_{v_0}=l_{v_0}$ for some $v_0<\footnote{we say $v_0< v_1$ if $v_0$ is a node before $v_1$ in the route and $v_0\neq v_1$.} v_1, v_0\in \mathcal{F}_i$.
        
      If not, we have 
       $\Delta_1 = min\big\{T_{v_1}-arr_{v_1}, \{l_{v}-T_v\}_{v< v_1, v\in \mathcal{F}_i}\big\}>0$.
      We can obtain a new feasible schedule $\mathcal{T}_1$ by delaying the service start times of all nodes in $\mathcal{F}_i$ that are before node $v_1$ (i.e., node $v$ such that $v< v_1, v\in \mathcal{F}_i$) to $T_v'= T_v + \Delta_1$. The excess user ride time of  $\mathcal{T}_1$ is at least $\Delta_1$ smaller than $\mathcal{T}$. It contradicts to our assumption that $\mathcal{T}$ is an optimal schedule (for example, see Example \ref{example1 of case 2});
        \item[(ii)] $T_{v_3}=e_{v_3}$ for some $v_3\geqslant v_2, v_3\in \mathcal{F}_i$.
        
         If not, we have $\Delta_2 = min\big\{T_{v_2}-arr_{v_2}, \{T_v-e_v\}_{v\geqslant v_2, v\in \mathcal{F}_i}\big\}>0$. We can obtain a new feasible schedule $\mathcal{T}_2$ by moving forward the service start times of all nodes in $\mathcal{F}_i$ that are after node $v_2$ (i.e., node $v$ such that $v\geqslant v_2, v\in \mathcal{F}_i$) to $T_v''= T_v - \Delta_2$. The excess user ride time of $\mathcal{T}_2$ is at least $\Delta_2$ smaller than $\mathcal{T}$. It contradicts to our assumption that $\mathcal{T}$ is an optimal schedule (for example, see Example \ref{example2 of case 2});
    \end{itemize}
    Based on (i) and (ii), assuming that $v_s, v_e$ are the first and the last node of $\mathcal{F}_i$, we derive that all the feasible schedules for $\mathcal{F}_i$ must satisfy the following two points:
    \begin{itemize}
        \item [(iii)] Since we have $arr_v= T_v$ for all $v<v_0 < v_1$ and $T_{v_0}=l_{v_0}$, any feasible schedules over $\mathcal{F}_i$ could not begin service at $v_s$ later than $T_{v_s}$ ($T_{v_s}$ is the latest possible service start time at $v_s$). Otherwise, it will surpass the latest time window $l_{v_0}$ at node $v_0$;
        \item [(iv)] Since we have $arr_v= T_v$ for all $v_2 \leqslant v_3 < v$ and $T_{v_3}=e_{v_3}$, any feasible schedules over $\mathcal{F}_i$ could not arrive at $v_e$ earlier than $arr_{v_e}$. 
    \end{itemize}
   
    Assuming that $\mathcal{T}_i^* = [\cdots,T_v^*,\cdots]_{v\in \mathcal{F}_i} $ is an optimal schedule of $\mathcal{F}_i$, and the arrival time at $v$ is $arr_v^* = T_{v-1}^* + t_{v-1,v}+s_{v-1}$. Now, we prove that the excess user ride time of $\mathcal{T}_i$ is the same as $\mathcal{T}_i^*$ using the above properties. Note that we are still under the condition that $arr_v< T_v$ for some $v\in \mathcal{F}_i$.
    
    According to (iii) and (iv), we have  $T_{v_s}^*\leqslant T_{v_s}, arr_{v_e}^*\geqslant arr_{v_e}$ for an optimal schedule $\mathcal{T}_i^*$ over $\mathcal{F}_i$. Clearly, $\mathcal{T}_i^*$ satisfies $ EU_{min}(\mathcal{T}_i^*)\leq EU_{min}(\mathcal{T}_i)$. Next, we will prove that $ EU_{min}(\mathcal{T}_i^*)= EU_{min}(\mathcal{T}_i)$. Then, we prove $\mathcal{T}_i$ is an optimal schedule over $\mathcal{F}_i$.
    Our proof contains two cases:
    \begin{enumerate}
        \item If $arr_v^* = T_v^*$ for all $v\in \mathcal{F}_i$:  As we have $T_{v_s}^*\leqslant T_{v_s}$, then $arr_{v_e}^*\leqslant arr_{v_e}$. Therefore,  we derive that $arr_{v_e}^*= arr_{v_e}$, $T_{v_s}^*= T_{v_s}$. As we assume in the condition that $arr_v^* = T_v^*$ for all $v\in \mathcal{F}_i$,  we must have $T_v = T_v^*$ for all $k\in \mathcal{F}_i$. It contradicts to our assumptions that $arr_v< T_v$ for some $v\in \mathcal{F}_i$. Therefore, this case will not happen;
        \item If $arr_v^* < T_v^*$ for some $v\in \mathcal{F}_i$: Then we can prove the same result as in (i) (ii) and (iii) for $T_v^*$ in the same manner. Then $T_{v_s}\leqslant T_{v_s}^*, arr_{v_e}\geqslant arr_{v_e}^*$ and thus we derive $T_{v_s}= T_{v_s}^*, arr_{v_e}= arr_{v_e}^*$.
        Then we have $EU_{min}(\mathcal{T}_i^*)=EU_{min}(\mathcal{T}_i)$.
       Otherwise, if $EU_{min}(\mathcal{T}_i^*)< EU_{min}(\mathcal{T}_i)$, we can obtain a new feasible schedule $\mathcal{T}'$ over $\mathcal{R}$ from $\mathcal{T}$ by replacing  $\mathcal{T}_i$ to $\mathcal{T}_i^*$, and $\mathcal{T}'$ has smaller excess user ride time than $\mathcal{T}$, which is a contradiction! 
    \end{enumerate}
\end{enumerate}
\end{proof}

For the sake of illustration, we take Example \ref{example1 of case 2} and \ref{example2 of case 2} to explain point (i) and (ii) of case 2, respectively.

\begin{exmp}\label{example1 of case 2}
Consider a fragment $\mathcal{F}=\{1+,2+,2-,1-\}$, the time window on each node, and the travel time for each arc is shown in Figure \ref{case2-1}. The direct travel time from node 1+ to node 1- is shown on the dashed line. Assume that the service time at each node is equal to zero and each request includes one passenger to be transported. We present two sets of service start times in the table, one is called schedule $\mathcal{B}$, and the other is called schedule $\mathcal{A}$. Also, we present the excess user ride times that are calculated from schedule $\mathcal{B}$ in the row named $R_{\mathcal{B}}$, and these are calculated from schedule $\mathcal{A}$ in the row named $R_{\mathcal{A}}$. We denote the service start time of schedule $\mathcal{B}$ at node $v$ as $\mathcal{B}_v$. Note that schedule $\mathcal{A}$ is an \textbf{optimal schedule}.

\begin{figure}[H]
\centering
\vspace{-4mm}
\includegraphics[width=12cm]{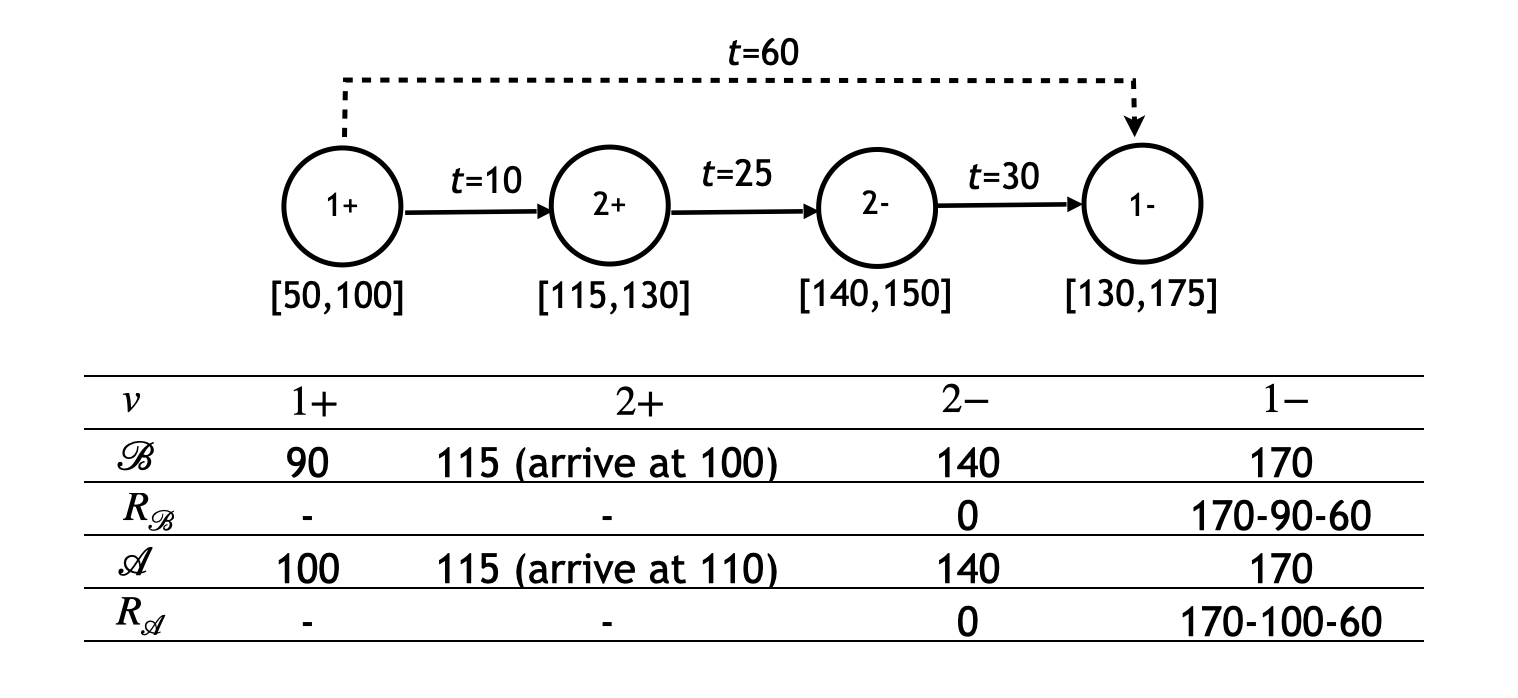}
\vspace{-6mm}
\caption{\centering Example of case 2 (i)}
\label{case2-1}
\end{figure}

In this example, node 2+ is the first node such that $arr_{2+} < \mathcal{A}_{2+}$ (i.e., node $v_1$ in the proof) and node 1+ is a node before node 2+ with $\mathcal{A}_{1+}=l_{1+}$ (i.e., node $v_0$ in the proof). Therefore, schedule $\mathcal{B}$ is a \textbf{conflicting schedule of case 2 (i)}. Schedule $\mathcal{B}$ is not excess-user-ride-time optimal, as one can further reduce excess user ride time of request 1 by delaying the service start time at node 1+ by $\Delta_1 = \min\{\mathcal{B}_{2+}-arr_{2+},l_{1+}-\mathcal{B}_{1+}\} = \min\{115-100,100-90\}$, as schedule $\mathcal{A}$ does.
\end{exmp}

\begin{exmp}[Example \ref{example1 of case 2} continued]\label{example2 of case 2}
We take the same problem settings as in Example \ref{example1 of case 2} to illustrate point (ii) of case 2. Note that schedule $\mathcal{A}$ is an \textbf{optimal schedule}.
    \begin{figure}[H]
    \centering
    \vspace{-4mm}
    \includegraphics[width=13cm]{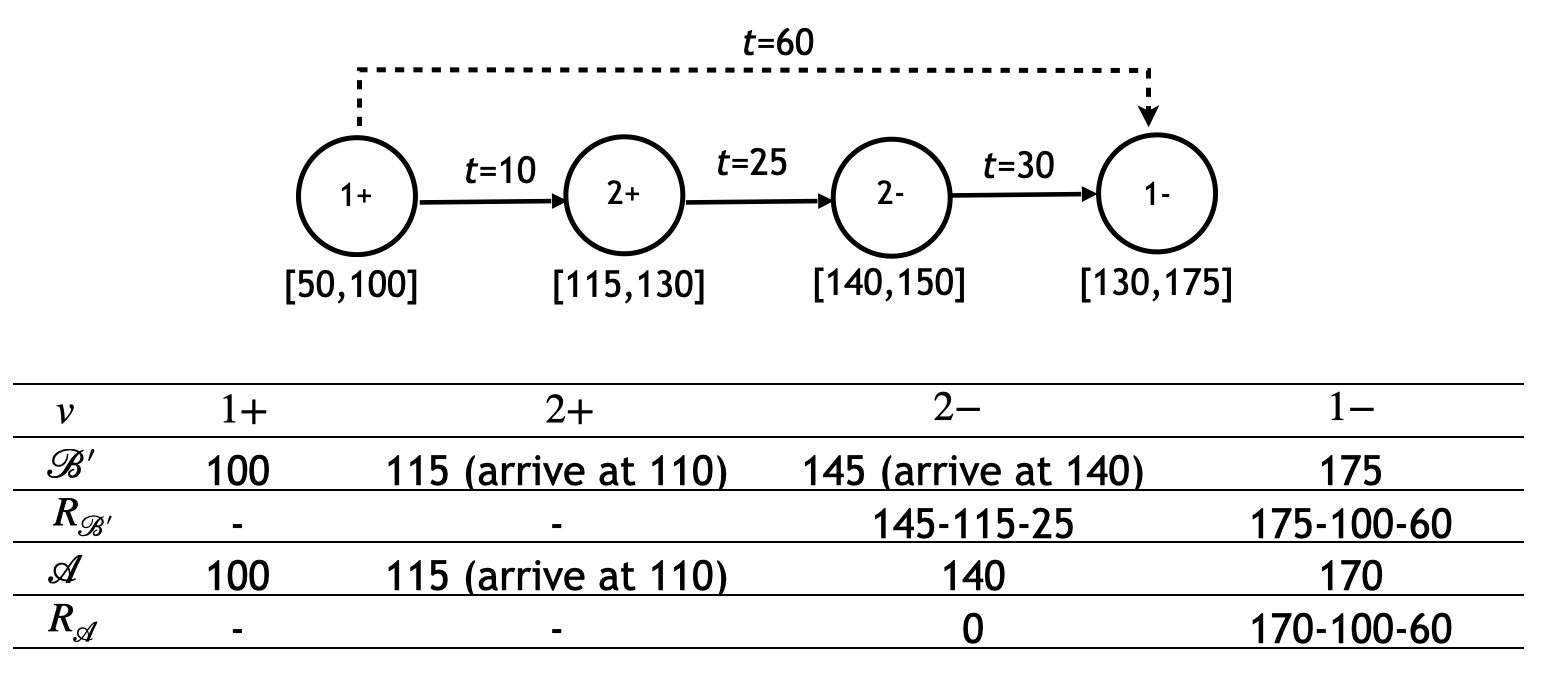}
    \vspace{-6mm}
    \caption{\centering Example of case 2 (ii)}
    \label{case2-2}
\end{figure}
In this example, node 2+ is the last node such that $arr_{2+} < \mathcal{A}_{2+}$ (i.e., node $v_2$ in the proof) and node 2- is a node after node 2+ with $\mathcal{A}_{2+} = e_{2+}$ (i.e., node $v_3$ in the proof). Therefore, schedule $\mathcal{B}'$ is a \textbf{conflicting schedule of case 2 (ii)}. Schedule $\mathcal{B}'$ is not an optimal schedule, as one can further reduce excess user ride times of request 1 and request 2 by moving forward the service start time at node 2- by: 
$$\Delta_2 = \min\{\mathcal{B}_{2+}'-arr_{2+},\{\mathcal{B}_{2-}'-e_{2-}, \mathcal{B}_{1-}'-e_{1-}\}\} = \min\{115-110,\{145-140,175-130\}\}.$$
\end{exmp}

\subsection{Linear Programming Model}
To calculate the minimum excess user ride time for a fragment contains two or more requests, the following Linear Programming (LP) model is invoked.
Let $P_\mathcal{F}$ denote the set of requests served on a fragment $\mathcal{F}$:

\begin{equation} \label{objective R}
    \quad  \min \sum\limits_{i \in P_\mathcal{F}}R_i
\end{equation}
s.t.

\begin{equation} \label{TW3}
    T_i + s_i + t_{i,j} \leqslant T_j, \quad \forall i \in \mathcal{F}, \quad idx_j = idx_i + 1, idx_i \neq |\mathcal{F}| 
\end{equation}
\begin{equation} \label{user ride 1}
    T_{n+i} - (T_i + s_i) \leqslant m_i, \quad \forall i \in P_\mathcal{F}
\end{equation}
\begin{equation} \label{user ride 2}
     T_{n+i} - T_i - s_i - t_{i,n+i} \leqslant R_i, \quad \forall i \in P_\mathcal{F}
\end{equation}
\begin{equation} \label{TW2}
     e_i \leqslant T_i \leqslant l_i, \quad \forall i \in \mathcal{F}
\end{equation}
\begin{equation}
      R_i \geqslant 0, \quad \forall i \in P_\mathcal{F}
\end{equation}

where $T_i$ denotes the service start time at node $i$, $idx_i$ is the index of node $i$ on the fragment. The objective function is to minimize the total excess user ride time of $\mathcal{F}$. Constraints (\ref{TW3}) are time window constraints. Constraints (\ref{user ride 1}) and constraints (\ref{user ride 2}) are user ride time constraints.

\section{Experimental results for fragment enumeration}
\label{preliminary for seg enumeration}
We provide all the details for enumerating all the fragments for each instance. In Table \ref{details seg}, $N_{frag}$ denotes the number of fragments generated, $Leg_{avg}$ and $Leg_{max}$ denote the average and maximum length of fragments, respectively. $N_{LP}$ represents the number of time LP is solved, and CPU is the total computational time for enumeration in seconds.

\begin{table}[!htp]
  \renewcommand\arraystretch{0.75}
  \centering
     \small
    \begin{threeparttable}
 \caption{Details of fragments enumeration for all the instances}
 \label{details seg}
   \begin{tabular}{c c c c c c}
        \hline
        & $N_{frag}$ & $Leg_{avg}$ & $Leg_{max}$ & $N_{LP}$ & CPU(s) \\
        \hline
        a2-16 &32  & 3.06 & 6  & 0 & 0.94\\
        a2-20 &51  & 3.41 & 6  & 1  & 0.23\\
        a2-24 &64  & 3.72 & 8  & 1  & 0.09\\
        a3-18 &71  & 4.25 & 8  & 4  & 0.04\\
        a3-24 &110 & 4.71 & 12 & 0  & 0.06\\
        a3-30 &89  & 3.66 & 8  & 0  & 0.12\\
        a3-36 &114 & 4.12 & 12 & 1  & 0.27\\
        a4-16 &78  & 4.51 & 8  & 4  & 0.04\\
        a4-24 &91  & 4.07 & 8  & 2 & 0.07\\
        a4-32 &206 & 5.58 & 12 & 3  & 0.20\\
        a4-40 &242 & 5.45 & 12 & 6 & 0.37\\
        a4-48 &355 & 5.33 & 12 & 15 & 0.61\\
        a5-40 &337 & 5.65 & 12 & 3 & 0.38\\
        a5-50 &659 & 8.25 & 24 & 33& 0.99\\
        \hline
        Avg   &178.5 &4.70 &10.57 &5.21 &0.32 \\
        \hline
        u2-16 &61   & 3.80 & 6  & 0  & 1.05\\
        u2-20 &180  & 5.26 & 12 & 7  & 0.32\\
        u2-24 &66   & 3.27 & 4  & 0  & 0.06\\
        u3-18 &78   & 3.95 & 8 & 0  & 0.04\\
        u3-24 &129  & 4.25 & 8  & 0  & 0.08\\
        u3-30 &255  & 5.06 & 8  & 19 & 0.29 \\
        u3-36 &276  & 5.14 & 12 & 12 & 0.30\\
        u4-16 &75  & 4.03 & 8 & 1  & 0.04\\
        u4-24 &57   & 3.19 & 6  & 0  & 0.05\\
        u4-32 &177  & 4.14 & 10 & 3  & 0.21\\
        u4-40 &149  & 4.01 & 8  & 2  & 0.26\\
        u4-48 &1177 & 9.01 & 18 & 7  & 1.69\\
        u5-40 &335  & 5.28 & 14 & 1  & 0.49\\
        u5-50 &584  & 6.13 & 14 & 6  & 0.96\\
        \hline
        Avg   &257.07 &4.75 &9.71 & 4.14 &0.42 \\
        \hline
        r5-60  &632   & 6.44  & 16 & 44  & 2.61 \\
        r6-48   & 4082  & 14.20 & 36 & 414  & 6.89\\
        r6-60   &809  & 6.58  & 18 & 40   & 1.65\\
        r6-72   &1080  & 7.12  & 22 & 36   & 2.51\\
        r7-56   &1089  & 7.92  & 18 & 83   & 1.70\\
        r7-70   &2340  & 8.32  & 18 & 183  & 4.14\\
        r7-84   &2892  & 11.66 & 30 & 405  & 7.77\\
        r8-64   &11694 & 18.23 & 42 & 3517 & 40.52\\
        r8-80   &5822  & 14.89 & 30 & 260  & 14.07\\
        r8-96   &3155  & 9.30  & 26 & 312  & 9.65\\
        \hline
        Avg   &3359.50 &10.47 &25.6 &526.4 &9.15 \\
        \hline
    \end{tabular}
    \end{threeparttable}
\end{table}

\section{Sensitivity Analysis and Parameter Tuning Results for $\theta_{max}$}
\label{detailed parameter tuning}

In this section, we present aggregated parameter tuning results for $\theta_{max}$ in Table \ref{t_max}. Then, we present the DA algorithm results with different settings of $\theta_{max}$ in Table \ref{detailed parameter tuning results}, where $F_{\theta_{max}}$ denotes the number of feasible solutions obtained for the associated instance among 10 runs of the DA algorithm with $\theta_{max}$.

\begin{table}[h!]
    \centering
     \small
    \begin{threeparttable}
    \caption{Sensitivity analysis for $\theta_{max}$ under different $\gamma$ cases on type-a instances}
    \label{t_max}
    \setlength\tabcolsep{7.5pt}
    \begin{tabular}{c c c c c c c c}
    \toprule
    $\theta_{max}$ & 0.6	&0.9 &1.2 &1.5	&1.8	&2.1	&2.4  \\
    \hline
    \textbf{$\gamma = 0.1$} \\
    \cline{1-1}
    $BC\%$ &0.11\%  & 0.10\%   & 0.19\%  & 0.32\%  & 0.29\%  & 0.28\%  & 0.51\%  \\
    $AC\%$ &0.49\%   & 0.53\%  & 0.74\%  & 0.83\%  & 0.82\%  & 0.94\%  & 1.07\%   \\
    $Q1\%$&0.23\%  & 0.30\%   & 0.43\%  & 0.54\%  & 0.56\%  & 0.66\%  & 0.81\%\\
    $Q3\%$&0.66\%  & 0.73\%  & 0.90\%   & 0.93\%  & 1.04\%  & 1.22\%  & 1.28\%  \\
    FeasRatio &140/140 &140/140 &140/140 &140/140 &140/140 &140/140 &140/140\\
    CPU (s) &83.93 & 77.43 & 78.52 & 80.09 & 81.16 & 82.12 & 83.42\\
    \hline
    \textbf{$\gamma = 0.4$} \\
    \cline{1-1}
    $BC\%$&0.19\%   & 0.27\%   & 0.27\%   & 0.40\%    & 0.49\%   & 0.70\%    & 0.63\% \\
    $AC\%$&NC    & 0.68\%   & 0.79\%   & 0.95\%   & 1.18\%   & 1.36\%   & 1.54\%\\
    $Q1\%$&0.31\%   & 0.49\%   & 0.57\%   & 0.65\%   & 0.84\%   & 0.96\%   & 1.11\%  \\
    $Q3\%$&0.72\%   & 0.84\%   & 0.97\%   & 1.21\%   & 1.5\%    & 1.68\%   & 1.83\%\\
     FeasRatio &139/140 &140/140 &140/140 &140/140 &140/140 &140/140 &140/140\\
    CPU (s) &121.34 & 116.97 & 119.03 & 121.72 & 122.97 & 125.65 & 127.80\\
    \hline
    \textbf{$\gamma = 0.7$} \\
    \cline{1-1}
    FeasRatio &85/140 &106/140 &106/140 &108/140 &112/140 &105/140 &106/140 \\
    CPU (s) &227.05 & 201.68 & 206.06 & 212.5 & 215.86 & 221.04 & 222.31\\
    \bottomrule
    \end{tabular}
  \end{threeparttable}
\end{table}

\begin{table}[!htp]
\renewcommand\arraystretch{0.8}
\setlength\tabcolsep{2.5pt}
 \caption{DA algorithm results on type-a instances with different settings of $\theta_{max}$}
 \label{detailed parameter tuning results}
 \vspace{-5mm}
    \begin{center}
    \footnotesize
    \begin{tabular}{c| c c|c c|c c|c c|c c|c c|c c}
    \hline
        &\multicolumn{2}{c}{$\theta_{max} = 0.6$} & \multicolumn{2}{c}{$\theta_{max} = 0.9$} &\multicolumn{2}{c}{ $\theta_{max} = 1.2$} &\multicolumn{2}{c}{ $\theta_{max} = 1.5$} & \multicolumn{2}{c}{$\theta_{max} = 1.8$} & \multicolumn{2}{c}{ $\theta_{max} = 2.1$} & \multicolumn{2}{c}{$\theta_{max} = 2.4$}\\
        \hline
        \textbf{$\gamma = 0.1$} & $BC_{0.6}\%$ & $AC_{0.6}\%$ & $BC_{0.9}\%$ & $AC_{0.9}\%$ & $BC_{1.2}\%$ & $AC_{1.2}\%$ & $BC_{1.5}\%$ & $AC_{1.5}\%$ & $BC_{1.8}\%$ & $AC_{1.8}\%$ & $BC_{2.1}\%$ & $AC_{2.1}\%$ & $BC_{2.4}\%$ & $AC_{2.4}\%$\\
        \hline
        a2-16 & 0    & 0    & 0    & 0    & 0    & 0    & 0    & 0    & 0    & 0    & 0    & 0    & 0    & 0 \\
        a2-20 &0    & 0    & 0    & 0    & 0    & 0    & 0    & 0    & 0    & 0    & 0    & 0    & 0    & 0\\
        a2-24 &0    & 0    & 0    & 0    & 0    & 0    & 0    & 0    & 0    & 0    & 0    & 0    & 0    & 0   \\
        a3-18 &0    & 0    & 0    & 0    & 0    & 0    & 0    & 0    & 0    & 0    & 0    & 0    & 0    & 0  \\
        a3-24 &0    & 0    & 0    & 0    & 0    & 0.74\% & 0    & 0.37\% & 0    & 0    & 0    & 0.15\% & 0    & 0.52\%\\
        a3-30 &0    & 0    & 0    & 0    & 0    & 0    & 0    & 0    & 0    & 0    & 0    & 0    & 0    & 0 \\
        a3-36 &0    & 0    & 0    & 0    & 0    & 0.12\% & 0    & 0    & 0    & 0.11\% & 0    & 0.10\%  & 0    & 0.62\%\\
        a4-16 &0    & 0    & 0    & 0    & 0    & 0    & 0    & 0    & 0    & 0    & 0    & 0    & 0    & 0\\
        a4-24 &0    & 0    & 0    & 0    & 0    & 0    & 0    & 0.70\%  & 0    & 0.67\% & 0    & 0.35\% & 0    & 0.70\%\\
        a4-32 &0    & 0.13\% & 0    & 1.04\% & 0.02\% & 0.79\% & 0.08\% & 1.31\% & 0    & 1.02\% & 0.06\% & 1.11\% & 0.09\% & 1.22\%\\
        a4-40 &0    & 0.62\% & 0    & 1.24\% & 0    & 1.23\% & 0.89\% & 1.38\% & 0    & 1.6\%  & 0    & 1.81\% & 1.24\% & 1.85\%\\
        a4-48 &0.34\% & 0.71\% & 0.30\%  & 0.94\% & 0.67\% & 1.78\% & 0.90\%  & 1.99\% & 0.82\% & 2.16\% & 1.17\% & 2.38\% & 1.55\% & 3.3\% \\
        a5-40 &0.44\% & 1.36\% & 0.29\% & 1.91\% & 0.42\% & 1.34\% & 0.78\% & 1.75\% & 0.96\% & 2.46\% & 1.16\% & 2.75\% & 1.16\% & 2.55\% \\
        a5-50 &0.71\% & 1.68\% & 0.74\% & 2.14\% & 1.61\% & 2.99\% & 1.83\% & 2.56\% & 2.33\% & 3.29\% & 1.57\% & 3.80\%  & 3.13\% & 3.98\% \\
         \hline
         Avg &0.11\% & 0.49\% & 0.10\%  & 0.53\% & 0.19\% & 0.74\% & 0.32\% & 0.83\% & 0.29\% & 0.82\% & 0.28\% & 0.94\% & 0.51\% & 1.07\% \\
        \hline
        \textbf{$\gamma = 0.4$} & $BC_{0.6}\%$ & $AC_{0.6}\%$ & $BC_{0.9}\%$ & $AC_{0.9}\%$ & $BC_{1.2}\%$ & $AC_{1.2}\%$ & $BC_{1.5}\%$ & $AC_{1.5}\%$ & $BC_{1.8}\%$ & $AC_{1.8}\%$ & $BC_{2.1}\%$ & $AC_{2.1}\%$ & $BC_{2.4}\%$ & $AC_{2.4}\%$\\
        \hline
        a2-16 &0     & 0     & 0     & 0     & 0     & 0     & 0     & 0     & 0     & 0    & 0     & 0     & 0     & 0\\
        a2-20 &0     & 0     & 0     & 0     & 0     & 0     & 0     & 0     & 0     & 0    & 0     & 0     & 0     & 0\\
        a2-24 &-0.29\% & -0.29\% & -0.29\% & -0.29\% & -0.29\% & -0.29\% & -0.29\% & -0.21\% & -0.29\% & 0.05\% & -0.29\% & -0.12\% & -0.29\% & 0.01\% \\
        a3-18 &0     & 0     & 0     & 0     & 0     & 0     & 0     & 0     & 0     & 0    & 0     & 0     & 0     & 0  \\
        a3-24 &0     & 0     & 0     & 0.27\%  & 0     & 0     & 0     & 0.23\%  & 0     & 0.74\% & 0     & 0.74\%  & 0     & 0.53\% \\
        a3-30 &-0.01\% & NC    & -0.01\% & -0.01\% & -0.01\% & -0.01\% & -0.01\% & -0.01\% & -0.01\% & 0.06\% & -0.01\% & 0.08\%  & -0.01\% & -0.01\%\\
        a3-36 &-0.22\% & -0.06\% & -0.16\% & 0.01\%  & -0.06\% & 0.05\%  & -0.03\% & 0.30\%   & -0.13\% & 0.94\% & -0.22\% & 0.64\%  & -0.06\% & 1.14\% \\
        a4-16 &0     & 0     & 0     & 0     & 0     & 0     & 0     & 0     & 0     & 0    & 0     & 0.12\% & 0     & 0.12\%  \\
        a4-24 &0     & 0.45\%  & 0     & 0.09\%  & 0.08\%  & 0.54\%  & 0.08\%  & 0.70\%   & 0     & 0.50\%  & 0.08\%  & 0.55\%  & 0     & 0.70\%\\
        a4-32 &0     & 0.20\%   & 0.22\%  & 1.26\%  & 0.22\%  & 0.75\%  & 0     & 1.0\%     & 0     & 1.64\% & 1.18\%  & 1.65\%  & 0.17\%  & 1.55\% \\
        a4-40 &0.27\%  & 0.98\%  & 0     & 1.28\%  & 0     & 1.41\%  & 0.27\%  & 1.70\%   & 0.99\%  & 1.61\% & 1.18\%  & 2.58\%  & 1.01\%  & 3.41\%\\
        a4-48 & 1.04\%  & 1.50\%   & 1.46\%  & 2.37\%  & 1.8\%   & 2.53\%  & 1.85\%  & 2.88\%  & 2.36\%  & 4.44\% & 2.87\%  & 4.02\%  & 3.17\%  & 4.54\%  \\
        a5-40 &0.35\%  & 1.56\%  & 0.55\%  & 1.19\%  & 0.86\%  & 1.40\%   & 0.55\%  & 1.52\%  & 1.21\%  & 2.24\% & 1.42\%  & 2.35\%  & 1.86\%  & 3.36\%\\
        a5-50 &1.46\%  & 2.03\%  & 2.07\%  & 2.84\%  & 1.19\%  & 3.45\%  & 3.14\%  & 4.30\%   & 2.71\%  & 3.94\% & 3.56\%  & 5.25\%  & 2.92\%  & 5.63\%\\
        \hline
        Avg &0.19\%  & NC   & 0.27\%  & 0.68\%  & 0.27\%  & 0.79\%  & 0.40\%   & 0.95\%  & 0.49\%  & 1.18\% & 0.7\%   & 1.36\%  & 0.63\%  & 1.54\% \\
        \hline
        \textbf{$\gamma = 0.7$} & $BC_{0.6}\%$ & $F_{0.6}$ & $BC_{0.9}\%$ & $F_{0.9}$ & $BC_{1.2}\%$ & $F_{1.2}$ & $BC_{1.5}\%$ & $F_{1.5}$ & $BC_{1.8}\%$ & $F_{1.8}$ & $BC_{2.1}\%$ & $F_{2.1}$ & $BC_{2.4}\%$ & $F_{2.4}$\\
        \hline
        a2-16 &0       & 9   & 0       & 10   & 0        & 10   & 0       & 10   & 0     & 10   & 0     & 10   & 0     & 10\\
        a2-20 &-       & 7   & -       & 8    & -        & 10    & -       & 7    & -     & 10   & -     & 7    & -     & 8\\
        a2-24 &0.85\%  & 1   & 0.85\%  & 8    & 0.85\%  & 8    & 0.85\%  & 8    & 0.82\%& 10   & 0.82\%  & 10   & 0.82\%  & 10  \\
        a3-18 &0       & 10  & 0       & 10   & 0        & 10   & 0       & 10   & 0     & 10   & 0     & 10   & 0     & 10 \\
        a3-24 &-0.63\% & 10  & -0.63\% & 10   & -0.33\%  & 10   & -0.63\% & 10   & -0.63\% & 10   & -0.33\% & 10   & -0.38\% & 10 \\
        a3-30 &-       & 5   & -       & 9    & -        & 6    & -     & 9    & -     & 10   & -     & 8    & -     & 8\\
        a3-36 &0.51\%  & 4   & 0       & 10   & 0.02\%        & 10   & 0.36\%  & 10   & 0.56\%  & 10   & 0.98\%  & 10   & 2.34\%  & 10 \\
        a4-16 &0       & 10  & 0       & 10   & 0        & 10   & 0     & 10   & 0     & 10   & 0     & 10   & 0     & 10\\
        a4-24 &0       & 10  & -0.49\% & 10   & -0.49\%        & 10   & 0     & 10   & -0.49\% & 10   & 0.03\%  & 10   & 0.03\%  & 10  \\
        a4-32 &-7.31\% & 10  & -6.40\% & 10   & -7.49\%  & 10   & -6.80\%  & 10   & -5.20\%  & 10   & -5.06\% & 10   & -6.07\% & 10 \\
        a4-40 &-       & 2   & -       & 1    & -        & 2    & -     & 4    & -     & 2    & NA   & 0    & NA   & 0\\
        a4-48 &NA      & 0   & NA      & 0    & NA       & 0    & NA   & 0    & NA   & 0    & NA   & 0    & NA   & 0 \\
        a5-40 &-4.43\% & 10  & -3.81\% & 10   & -5.23\%  & 10   & -2.44\% & 10   & -1.25\% & 10   & -3.00\%    & 10   & -1.47\% & 10  \\
        a5-50 &NA      & 0   & NA      & 0    & NA       & 0    & NA   & 0    & NA   & 0    & NA   & 0    & NA   & 0 \\
        \hline
    \end{tabular}
    \vspace{-2mm}
    \begin{tablenotes}
    \footnotesize
     \item ``-" indicates we obtain new best solution on previously unsolved instance and the gap cannot be calculated.
  \end{tablenotes}
    \end{center}
\end{table}

\section{Sensitivity analysis of increasing the maximum number of charging visits per station}
\label{unlimited visits}
Tables \ref{mul cordeau}, \ref{mul uber}, and \ref{mul ropke} present the DA algorithm results on type-a, type-u, and type-r instances with allowing at-most-one, at-most-two, at-most-three, and unlimited visits per recharging station, respectively. The maximum allowed charging visits per station is denoted by $n_{as}$.
 
\begin{table}[!htp]
\setlength\tabcolsep{3.5pt}
 \caption{Solution quality and performance on type-a instances when increasing the maximum number of charging visits per station}
 \label{mul cordeau}
    \begin{center}
    \footnotesize
    \begin{tabular}{c c c c| c c c| c c c| c c c c c}
    \hline
        &\multicolumn{3}{c}{DA with $n_{as} = 1$} & \multicolumn{3}{c}{DA with $n_{as} = 2$} &\multicolumn{3}{c}{DA with $n_{as} = 3$} &\multicolumn{5}{c}{DA with $n_{as} = \infty$}\\
        \hline
        \textbf{$\gamma = 0.4$} & $BC$ & $AC$ &CPU & $BC_2$ & $AC_2$ &CPU$_2$ & $BC_3$ & $AC_3$ &CPU$_3$ & $BC_{\infty}$ & $AC_{\infty}$ &CPU$_{\infty}$ &$N_{max}^s$ &$\overline{N_{avg}^s}$\\
        \hline
        a2-16 &\textbf{237.38} & 237.38 & 52.85  & \textbf{237.38} & 237.38 & 52.65  & \textbf{237.38} & 237.38 & 53.07  & \textbf{237.38} & 237.38 & 50.65 &1 &1.0\\
        a2-20 &\textbf{280.70} & 280.70 & 140.70 & \textbf{280.70} & 280.70 & 148.12 & \textbf{280.70} & 280.70 & 141.97 &\textbf{280.70} & 280.70 & 144.92 &1 &1.0\\
        a2-24 &347.04 & 347.04 & 230.99 & \textbf{346.28} & 346.28 & 286.96 & \textbf{346.28} & 346.28 & 284.47 & \textbf{346.28} & 346.28 & 265.80 &2 &3.0\\
        a3-18 &\textbf{236.82} & 236.82 & 26.30  & \textbf{236.82} & 236.82 & 26.93  & \textbf{236.82} & 236.82 & 26.43  & \textbf{236.82} & 236.82 & 25.36 &0 &0.0\\
        a3-24 &\textbf{274.80} & 274.80 & 67.85  & \textbf{274.80} & 274.80 & 71.07  & \textbf{274.80} & 274.80 & 69.48  & \textbf{274.80} & 274.80 & 66.66  &1 &0.8\\
        a3-30 &\textbf{413.34} & 413.34 & 88.67  & \textbf{413.34} & 413.34 & 104.70 & \textbf{413.34} & 413.34 & 106.13 & \textbf{413.34} & 413.34 & 103.54 &1 &2.0\\
        a3-36 &483.06 & 483.86 & 157.79 & \textbf{481.17} & 481.46 & 255.25 & \textbf{481.17} & 481.17 & 264.23 & \textbf{481.17} & 481.17 & 248.95 &3 &3.0\\
        a4-16 &\textbf{222.49} & 222.49 & 19.39  & \textbf{222.49} & 222.49 & 19.71  & \textbf{222.49} & 222.49 & 19.08  & \textbf{222.49} & 222.49 & 17.78  &0 &0.0 \\
        a4-24 &\textbf{311.03} & 311.65 & 31.97  & \textbf{311.03} & 311.65 & 31.54  & \textbf{311.03} & 311.65 & 31.15  & \textbf{311.03} & 311.65 & 29.53 &0 &0.0 \\
        a4-32 &\textbf{394.26} & 397.21 & 62.95  & \textbf{394.26} & 397.31 & 65.66  & \textbf{394.26} & 397.21 & 63.85  & \textbf{394.26} & 397.27 & 61.71 &1 &1.0 \\
        a4-40 &\textbf{453.84} & 459.46 & 116.65 & \textbf{453.84} & 459.18 & 125.28 & \textbf{453.84} & 459.11 & 116.86 & \textbf{453.84} & 458.74 & 121.04 &1 &0.7 \\
        a4-48 &558.11 & 563.47 & 177.51 & 558.18 & 564.63 & 235.32 & \textbf{557.86} & 564.21 & 238.60 & 558.96 & 564.86 & 231.45  &2 &2.7 \\
        a5-40 &416.25 & 420.32 & 72.64  & 415.62 & 420.09 & 71.75  & \textbf{415.43} & 420.16 & 72.01  & 415.79 & 419.82 & 70.78 &0 &0.1 \\
        a5-50 &567.54 & 574.56 & 162.82 & \textbf{564.90} & 575.04 & 190.93 & 567.40 & 574.64 & 189.18 & 567.13 & 574.28 & 184.43  &1 &1.6 \\
        \hline
        Avg & &  &100.65   &  &  &120.42 & &  &119.75   &  &  &115.90 &1.0 &1.2 \\
        \hline
        \textbf{$\gamma = 0.7$} & $BC$ & $AC$ &CPU & $BC_2$ & $AC_2$ &CPU$_2$ & $BC_3$ & $AC_3$ &CPU$_3$ & $BC_{\infty}$ & $AC_{\infty}$ &CPU$_{\infty}$  &$N_{max}^s$ &$\overline{N_{avg}^s}$\\
        \hline
        a2-16 &\textbf{240.66} & 240.66 & 95.75  & \textbf{240.66} & 240.66 & 125.10 & \textbf{240.66} & 240.66 & 124.90 & \textbf{240.66} & 240.66 & 119.29 & 2  &3.0 \\
        a2-20 &293.27 & 294.11 & 172.77 & 286.52 & 286.52 & 331.90 & \textbf{285.86} & 285.86 & 327.01 & 286.52 & 288.89 & 316.22  & 2  &3.6  \\
        a2-24 &353.18 &NA & 206.58 & 352.25 & 363.17 & 373.77 & \textbf{350.49} & 361.02 & 390.86 & 354.38 & 374.68 & 357.33 &2  &3.9  \\
        a3-18 &240.58 & 240.58 & 58.30  & \textbf{238.82} & 238.82 & 70.27  & \textbf{238.82} & 238.82 & 69.52  & \textbf{238.82} & 238.82 & 65.89  & 3  &4.0  \\
        a3-24 &275.97 & 277.43 & 123.71 & \textbf{275.20} & 275.20 & 154.90 & \textbf{275.20} & 275.94 & 155.39 & \textbf{275.20} & 275.20 & 150.02  & 2  &2.9  \\
        a3-30 &424.93 & 436.20 & 77.73  & 416.87 & 417.90 & 173.80 & \textbf{415.71} & 417.35 & 176.38 & \textbf{415.71} & 417.07 & 170.95  & 3 &4.7  \\
        a3-36 &494.04 & 502.27 & 125.42 & 486.36 & 487.34 & 332.47 & \textbf{484.85} & 487.59 & 350.73 & \textbf{484.85} & 487.91 & 343.02  & 3  &4.8  \\
        a4-16 &223.13 & 223.13 & 31.32  & \textbf{222.49} & 223.13 & 33.40  & \textbf{222.49} & 222.49 & 36.24  &\textbf{222.49} & 222.49 & 31.37  &2  &1.7  \\
        a4-24 &316.65 & 318.31 & 53.73  & \textbf{315.98} & 317.99 & 74.82  & \textbf{315.98} & 317.99 & 80.77  & \textbf{315.98} & 317.99 & 70.97  &2  &2.7 \\
        a4-32 &397.87 & 405.85 & 71.44  & 395.84 & 402.85 & 127.78 & \textbf{394.99} & 402.38 & 142.98 & \textbf{394.94} & 401.82 & 123.77 & 4 &3.7 \\
        a4-40 &479.02 & NA & 114.74 & 458.98 & 467.15 & 235.88 & 458.73 & 465.04 & 250.11 & \textbf{458.52} & 467.60 & 226.05  & 3 &4.6 \\
        a4-48 &582.22 & NA & 164.39 & 569.23 & 576.26 & 379.04 & \textbf{566.26} & 577.30 & 434.97 & 568.08 & 575.96 & 403.27  & 2 &5.3 \\
        a5-40 &424.26 & 436.94 & 97.51  & 417.35 & 424.29 & 153.00 & \textbf{416.89} & 423.96 & 169.49 & 419.33 & 425.29 & 149.77  & 4 &4.3 \\
        a5-50 &603.24 & NA & 158.39 & 583.37 & 590.81 & 320.55 & \textbf{576.54} & 589.38 & 367.00 & 579.15 & 588.98 & 352.73  & 4 &5.7 \\
        \hline
        Avg & &  &110.84   &  &  &206.19 & &  &219.74   &  &  &205.76  &2.7 &3.9 \\
        \hline
    \end{tabular}
    \end{center}
\end{table}

\begin{table}[!htp]
 \renewcommand\arraystretch{0.8}
 \setlength\tabcolsep{1pt}
 \caption{Solution quality and performance on type-u instances when increasing the maximum number of charging visits per station}
 \vspace{-2mm}
 \label{mul uber}
    \begin{center}
    \footnotesize
    \begin{tabular}{c c c c c|c c c c|c c c c|c c c c c}
    \hline
        &\multicolumn{4}{c}{DA with $n_{as} = 1$} & \multicolumn{4}{c}{DA with $n_{as} = 2$} &\multicolumn{4}{c}{DA with $n_{as} = 3$} &\multicolumn{5}{c}{DA with $n_{as} = \infty$}\\
        \hline
        \textbf{$\gamma = 0.1$} & $BC$ & $AC$ &CPU & $BC'$ & $BC_2$ & $AC_2$ &CPU$_2$ & $BC_2'$ & $BC_3$ & $AC_3$ &CPU$_3$ & $BC_3'$ & $BC_{\infty}$ & $AC_{\infty}$ &CPU$_{\infty}$ &$N_{max}^s$ &$\overline{N_{avg}^s}$\\
        \hline
        u2-16 &\textbf{57.61}  & 57.61  & 120.06 &57.61$^*$ & \textbf{57.61}  & 57.61  & 124.25 &57.61$^*$ &  \textbf{57.61}  & 57.61  & 126.90 &57.61$^*$ & \textbf{57.61}  & 57.61  & 182.64   & 1 &0.8\\
        u2-20 &\textbf{55.59}  & 56.34  & 401.82 &55.59$^*$ & \textbf{55.59}  & 55.59  & 421.77 &55.59$^*$ & \textbf{55.59}  & 55.59  & 440.65  &55.59$^*$ & \textbf{55.59}  & 55.59  & 642.95  & 1 &0.6\\
        u2-24 &\textbf{90.73$^*$}  & 90.84  & 599.73 &\textit{91.27$^*$} & \textbf{90.73$^*$}  & 90.73  & 572.00 &\textit{91.27$^*$} & \textbf{90.73$^*$}  & 90.73  & 592.75  &\textit{91.27$^*$} & \textbf{90.73}  & 90.73  & 1021.42  &  1 &2.9\\
        u3-18 &\textbf{50.74}  & 50.74  & 108.32 &50.74$^*$ & \textbf{50.74}  & 50.74  & 111.63 &50.74$^*$ & \textbf{50.74}  & 50.74  & 112.69  &50.74$^*$ & \textbf{50.74}  & 50.74  & 172.79  &  0 &0.0\\
        u3-24 &\textbf{67.56}  & 68.16  & 111.49 &67.56$^*$ & \textbf{67.56}  & 68.16  & 115.43 &67.56$^*$ & \textbf{67.56}  & 68.16  & 117.24  &67.56$^*$ & \textbf{67.56}  & 68.16  & 173.10  &  0 &0.0\\
        u3-30 &\textbf{76.75}  & 77.80  & 174.11 &76.75$^*$ & \textbf{76.75}  & 77.55  & 182.98 &76.75$^*$ & \textbf{76.75}  & 77.55  & 168.08  &76.75$^*$ & \textbf{76.75}  & 77.55  & 268.96  &  0 &0.3\\
        u3-36 &104.27 & 105.42 & 420.72 &104.04$^*$ & 104.27 & 105.45 & 578.30 &104.04$^*$ & 104.27 & 106.10 & 552.64  &104.04$^*$ & 104.27 & 105.48 & 775.41  & 1 &1.7\\
        u4-16 &\textbf{53.58}  & 53.58  & 51.37  &53.58$^*$ & \textbf{53.58}  & 53.58  & 51.14  &53.58$^*$ & \textbf{53.58}  & 53.58  & 49.18   &53.58$^*$ & \textbf{53.58}  & 53.58  & 72.84 & 0 &0.0\\
        u4-24 &90.13  & 90.85  & 55.26  &89.83$^*$ & 89.91  & 90.85  & 57.23  &89.83$^*$ & 90.08  & 90.85  & 56.87 &89.83$^*$  & 90.08  & 90.85  & 79.82 & 1 &0.5\\
        u4-32 &\textbf{99.29}  & 99.42  & 119.12 &99.29$^*$ & \textbf{99.29}  & 99.42  & 114.88 &99.29$^*$ & \textbf{99.29}  & 99.42  & 118.06  &99.29$^*$ & \textbf{99.29}  & 99.42  & 162.58 & 0 &0.4\\
        u4-40 &\textbf{133.11} & 135.18 & 154.00 &133.11$^*$ & \textbf{133.11} & 135.34 & 163.78 &133.11$^*$ & 133.14 & 135.21 & 159.92  &133.11$^*$ & \textbf{133.11} & 135.23 & 216.58 & 1 &1.7 \\
        u4-48 &\textbf{147.75$^*$} & 149.69 & 840.96 &148.30 & \textbf{147.73$^*$} & 149.89 & 917.71 &148.37 & \textbf{147.43$^*$} & 149.52 & 902.77 & 149.14 & \textbf{147.33$^{**}$} & 149.37 & 1403.39 & 2 &2.9\\
        u5-40 &\textbf{121.86} & 123.38 & 113.81 &121.86 & \textbf{121.86} & 123.54 & 116.57 &121.86 & \textbf{121.86} & 123.74 & 118.30 &121.86 & \textbf{121.86} & 123.59 & 149.98 & 1 &0.8\\
        u5-50 &144.22 & 145.63 & 245.52 &143.10 & 143.27 & 145.73 & 258.43 &142.83 & 143.51 & 145.91 & 279.38 &142.83 & \textbf{143.14$^{**}$} & 146.05 & 393.68 & 1 &1.5\\
        \hline
        Avg & &  &251.16 &   &  &  &270.43 & & &  &271.10 &  &  &  &408.30 &0.7 &1.0\\
        \hline
        \textbf{$\gamma = 0.4$} & $BC$ & $AC$ &CPU & $BC'$ & $BC_2$ & $AC_2$ &CPU$_2$ & $BC_2'$ & $BC_3$ & $AC_3$ &CPU$_3$ & $BC_3'$ & $BC_{\infty}$ & $AC_{\infty}$ &CPU$_{\infty}$ &$N_{max}^s$ &$\overline{N_{avg}^s}$\\
        \hline
        u2-16 &\textbf{57.65}  & 57.65  & 156.61  &57.65$^*$ & \textbf{57.65}  & 57.65  & 171.29  &57.65$^*$ & \textbf{57.65}  & 57.65  & 168.11 &57.65$^*$ & \textbf{57.65}  & 57.65  & 276.29 & 1 &2.0\\
        u2-20 &\textbf{56.34}  & 56.34  & 606.64  &56.34$^*$ & \textbf{56.34}  & 56.34  & 690.19  &56.34$^*$ & \textbf{56.34}  & 56.34  & 682.00 &56.34$^*$ & \textbf{56.34}  & 56.34  & 1006.29 &1 &2.0 \\
        u2-24 &\textbf{91.24$^*$}  & 91.72  & 817.79  &\textit{91.63$^*$} & \textbf{91.14$^*$}  & 91.43  & 836.02  &\textit{91.27$^*$} & \textbf{91.14$^*$}  & 91.43  & 885.85 &\textit{91.27$^*$} & 91.16  & 91.17  & 1399.38 & 2 &3.3\\
        u3-18 &\textbf{50.74}  & 50.74  & 124.95  &50.74$^*$ & \textbf{50.74}  & 50.74  & 129.61  &50.74$^*$ & \textbf{50.74}  & 50.74  & 133.92 &50.74$^*$ & \textbf{50.74}  & 50.74  & 213.60 & 1 &1.1\\
        u3-24 &\textbf{67.56}  & 68.16  & 141.01  &67.56$^*$ & 67.86  & 68.06  & 145.32  &67.56$^*$ & 67.67  & 68.16  & 153.68 &67.56$^*$ & \textbf{67.56} & 68.16  & 214.32  & 1 &1.2 \\
        u3-30 &\textbf{76.75}  & 77.93  & 285.81  &76.75$^*$ & \textbf{76.75}  & 78.13  & 306.82  &76.75$^*$ & \textbf{76.75}  & 78.28  & 298.28 &76.75$^*$ & \textbf{76.75}  & 77.85  & 420.10 & 1 &2.1\\
        u3-36 &104.49 & 106.37 & 898.90  &104.06$^*$ & \textbf{104.06} & 106.68 & 1038.76 &104.06$^*$ & 104.69 & 106.57 & 1078.92&104.06$^*$ & 104.31 & 106.07 & 1589.46 & 2 &3.5\\
        u4-16 &\textbf{53.58}  & 53.58  & 60.52   &53.58$^*$ & \textbf{53.58}  & 53.58  & 62.49   &53.58$^*$ & \textbf{53.58}  & 53.58  & 63.21  &53.58$^*$ & \textbf{53.58}  & 53.58  & 85.00 & 0 &0.0\\
        u4-24 &90.72  & 91.00  & 65.57   &89.83$^*$ & 90.21  & 90.90  & 68.48   &89.83$^*$ & 90.13  & 90.90  & 70.04  &89.83$^*$ & \textbf{90.08$^{**}$}  & 90.85  & 91.67 & 2 &2.1\\
        u4-32 &\textbf{99.29}  & 99.42  & 156.27  &99.29$^*$ & \textbf{99.29}  & 99.42  & 166.76  &99.29$^*$ & \textbf{99.29}  & 99.42  & 162.08 &99.29$^*$ &\textbf{99.29}  & 99.42  & 230.20 & 1 &2.7\\
        u4-40 &\textbf{133.78$^*$} & 135.83 & 303.06  &\textit{133.91$^*$} & \textbf{133.61$^*$} & 135.75 & 318.07  &\textit{133.68$^*$} & 134.23 & 136.16 & 326.55 &134.01 & \textbf{133.36$^{**}$} & 136.19 & 457.33 & 2 &4.2 \\
        u4-48 &\textbf{148.48$^*$} & 150.81 & 1390.74 &NA & \textbf{148.18$^*$} & 150.53 & 1247.04 &150.96 & \textbf{148.23$^*$} & 150.21 & 1454.38&150.78 & \textbf{147.75$^{**}$} & 149.71 & 2050.93 & 2 &4.9 \\
        u5-40 &\textbf{121.96$^*$} & 123.63 & 160.80  &122.23 & \textbf{121.96$^*$} & 123.50 & 163.39  &122.22 & \textbf{121.96} & 123.77 & 166.90 &121.96 & \textbf{121.96} & 123.94 & 237.16 & 1 &3.3 \\
        u5-50 &143.68 & 146.60 & 391.46  &143.14 & 143.78 & 146.36 & 401.78  &142.83 & 143.50 & 146.21 & 415.65 &143.48 & \textbf{143.42$^{**}$} & 145.65 & 619.05 & 1 &4.1 \\
        \hline
        Avg & &  &397.15 &  &  &  &410.43 &  & &  &432.83 &  &  &  &835.06 &1.3  &2.6 \\
        \hline
        \textbf{$\gamma = 0.7$} & $BC$ & $AC$ &CPU & $BC'$ & $BC_2$ & $AC_2$ &CPU$_2$ & $BC_2'$ & $BC_3$ & $AC_3$ &CPU$_3$ & $BC_3'$ & $BC_{\infty}$ & $AC_{\infty}$ &CPU$_{\infty}$ &$N_{max}^s$ &$\overline{N_{avg}^s}$\\
        \hline
        u2-16 &\textbf{59.19}  & 60.01  & 419.57  &59.19$^*$ & \textbf{58.17}  & 58.17  & 460.44  &58.17$^*$ &\textbf{58.17}  &58.17   &530.24  &58.17$^*$ & 58.75  & 59.46  & 663.32 & 2 &3.3 \\
        u2-20 &\textbf{56.86}  & 58.39  & 1527.60 &56.86$^*$ & \textbf{56.86}  & 58.03  & 1561.63 &56.86$^*$ & \textbf{56.86}  & 57.98  & 1583.70 &56.86$^*$ & \textbf{56.86}  & 58.39  & 2619.96 & 1 &2.8 \\
        u2-24 &\textbf{92.84$^*$}  & 99.38  & 1065.06 &NA & \textbf{92.43$^*$}  & 105.67 & 1307.99 &97.50 & \textbf{92.43$^*$}  & 101.95 & 1529.29 &NA & 92.77  & 100.36 & 2090.28 & 2 &5.0 \\
        u3-18 &\textbf{50.99}  & 50.99  & 206.92  &50.99$^*$ & \textbf{50.99}  & 50.99  & 206.48  &50.99$^*$ & \textbf{50.99}  & 50.99  & 217.78  &50.99$^*$ & \textbf{50.99}  & 50.99  & 301.43  & 1 &3.0 \\
        u3-24 &\textbf{68.39}  & 68.44  & 375.75  &68.39$^*$ & 68.24  & 68.39  & 389.47  &68.06$^*$ & 68.24  & 68.51  & 419.27  &68.06$^*$ & \textbf{68.06$^{**}$}  & 68.41  & 544.52  &  2 &3.8 \\
        u3-30 &\textbf{77.94$^*$}  & 79.37  & 1094.81 &\textit{78.14$^*$} & \textbf{77.94$^*$}  & 79.09  & 1132.97 &78.16 & \textbf{77.94$^*$}  & 79.02  & 1293.92 &78.16 & \textbf{77.83$^{**}$}  & 79.11  & 1595.22  & 2 &4.2 \\
        u3-36 &106.00 & 107.57 & 1606.43 &105.79 & \textbf{106.39$^*$} & 107.62 & 1521.37 &107.65 & 106.39 & 107.07 & 1605.03 &106.18 & \textbf{105.98$^{**}$} & 106.95 & 2690.77  &  2 &4.5 \\
        u4-16 &\textbf{53.87}  & 53.87  & 96.90   &53.87$^*$ & \textbf{53.87}  & 53.87  & 100.33  &53.87$^*$ & \textbf{53.87}  & 53.87  & 103.29  &53.87$^*$ & \textbf{53.87}  & 53.87  & 133.65 &  1 & 3.0\\
        u4-24 &90.07  & 90.97  & 254.45  &89.96$^*$ & 89.96  & 90.97  & 263.46  &89.83$^*$ & 89.91  & 90.97  & 282.62  &89.83 & \textbf{89.83$^{**}$}  & 90.72  & 375.00 &  2 &3.9 \\
        u4-32 &\textbf{99.50}  & 101.09 & 325.31  &99.50$^*$ & \textbf{99.50}  & 99.95  & 321.35  &99.50$^*$ & \textbf{99.50}  & 100.34 & 342.44  &99.50$^*$ & \textbf{99.50}  & 100.28 & 526.67 & 1 &4.6 \\
        u4-40 &\textbf{136.08$^*$} & 138.98 & 708.04  &NA & \textbf{134.98$^*$} & 138.37 & 731.95  &137.49 & \textbf{135.38$^*$} & 138.01 & 730.23  &137.61 & \textbf{134.94$^{**}$} & 136.20 & 971.29 & 2 &5.5 \\
        u4-48 &\textbf{152.58$^*$} & 162.62 & 1958.80 &NA & \textbf{150.55$^*$} & 154.19 & 1962.85 &NA & \textbf{151.57$^*$} & 155.36 & 1955.60 &NA & \textbf{149.51$^{**}$} & 152.90 & 2907.41 & 3 &6.3 \\
        u5-40 &\textbf{123.52$^*$} & 126.10 & 359.59  &NA & \textbf{124.04$^*$} & 126.08 & 385.25  &125.14 & \textbf{123.71$^*$} & 125.63 & 401.18  &124.18 & \textbf{123.32$^{**}$} & 125.15 & 506.11 & 2 &5.4 \\
        u5-50 &\textbf{143.51$^*$} & 149.52 & 922.19  &144.36 & \textbf{144.24$^*$} & 148.13 & 923.51  &164.19 & \textbf{143.51$^*$} & 148.53 & 1001.25 &144.10 & \textbf{142.89$^{**}$} & 146.10 & 1165.39 & 2 & 6.1\\
        \hline
        Avg & &  &780.10 &   &  &  &804.93 & & &  &856.84  &  &  &  &1220.79 & 1.8 &4.4 \\
        \hline
    \end{tabular}
    \end{center}
\end{table}

\begin{table}[!htp]
\setlength\tabcolsep{2.5pt}
 \caption{Solution quality and performance on type-r instances when increasing the maximum number of charging visits per station}
 \label{mul ropke}
    \begin{center}
    \footnotesize
   \begin{tabular}{c c c c| c c c| c c c| c c c c c}
    \hline
        &\multicolumn{3}{c}{DA with $n_{as} = 1$} & \multicolumn{3}{c}{DA with $n_{as} = 2$} &\multicolumn{3}{c}{DA with $n_{as} = 3$} &\multicolumn{5}{c}{DA with $n_{as} = \infty$}\\
        \hline
        \textbf{$\gamma = 0.1$} & $BC$ & $AC$ &CPU & $BC_2$ & $AC_2$ &CPU$_2$ & $BC_3$ & $AC_3$ &CPU$_3$ & $BC_{\infty}$ & $AC_{\infty}$ &CPU$_{\infty}$ &$N_{max}^s$ &$\overline{N_{avg}^s}$\\
        \hline
        r5-60  &691.83  & 706.20  & 178.44 & 689.75  & 703.86  & 175.42 & 688.52  & 706.91  & 180.86 & \textbf{687.68}  & 705.59  & 171.75  &0 &0.2\\
        r6-48  &506.72  & 512.69  & 229.31 & \textbf{506.45}  & 513.62  & 241.23 & 507.03  & 513.63  & 231.32 & 506.91  & 514.15  & 241.89  &0 &0.0\\
        r6-60  &692.00  & 700.15  & 127.03 & \textbf{690.15}  & 701.15  & 133.74 & 692.24  & 701.86  & 137.18 & 691.07  & 702.09  & 128.33  &0 &0.0\\
        r6-72  &777.44  & 794.69  & 208.39 & 776.68  & 795.41  & 212.78 & \textbf{775.93}  & 793.96  & 208.77 & 777.46  & 795.14  & 210.51 & 1 &0.1\\
        r7-56  &\textbf{613.10}  & 624.51  & 88.20  & 614.61  & 623.65  & 91.27  & 615.61  & 623.52  & 84.50  & 614.18  & 622.69  & 87.32  & 0 &0.0\\
        r7-70  &760.90  & 778.84  & 209.76 & 761.16  & 776.92  & 212.08 & 761.25  & 778.05  & 202.26 & \textbf{760.10}  & 777.10  & 202.03  & 0 &0.0\\
        r7-84  &889.38  & 904.88  & 322.66 & \textbf{884.43}  & 903.96  & 318.05 & 890.47  & 905.78  & 339.95 & 885.89  & 905.13  & 300.21  & 0 &0.1\\
        r8-64  &641.99  & 652.59  & 612.06 & \textbf{640.05}  & 653.65  & 645.07 & 642.09  & 653.44  & 773.82 & 640.24  & 653.81  & 647.97  &0 &0.0 \\
        r8-80  &803.52  & 828.67  & 357.75 & 807.04  & 826.91  & 366.82 & \textbf{799.00}  & 826.71  & 376.87 & 804.02  & 826.92  & 372.21  & 0 &0.0\\
        r8-96  &1053.11 & 1080.80 & 363.46 & 1052.19 & 1078.29 & 358.23 & 1064.64 & 1081.49 & 377.77 & \textbf{1049.98} & 1077.21 & 366.73  &0 &0.4 \\
        \hline
        Avg  & &  &269.71   &  &  &275.47 & &  &291.33   &  &  &272.90 &0.1 &0.1 \\
        \hline
        \textbf{$\gamma = 0.4$} & $BC$ & $AC$ &CPU & $BC_2$ & $AC_2$ &CPU$_2$ & $BC_3$ & $AC_3$ &CPU$_3$ & $BC_{\infty}$ & $AC_{\infty}$ &CPU$_{\infty}$ &$N_{max}^s$ &$\overline{N_{avg}^s}$\\
        \hline
        r5-60  &697.97 & 718.44 & 293.25 & 703.00  & 721.56  & 308.94 & 692.84  & 710.40  & 288.01 & \textbf{691.72}  & 709.78  & 285.00  & 2 &3.0\\
        r6-48  &506.91 & 514.46 & 257.59 & \textbf{506.45}  & 511.62  & 248.38 & 506.75  & 511.00  & 258.81 & 507.25  & 514.64  & 255.83  & 0 &0.1\\
        r6-60  &694.78 & 706.07 & 173.43 & 693.80  & 706.11  & 175.96 & 693.03  & 703.13  & 174.80 & \textbf{692.83}  & 701.86  & 174.24  & 1 &1.7 \\
        r6-72  &799.60 & 821.17 & 349.98 & 795.88  & 814.03  & 342.96 & \textbf{776.17}  & 800.29  & 336.47 & 781.22  & 801.86  & 342.33  & 1 &3.3 \\
        r7-56  &613.66 & 624.40 & 99.91  & \textbf{612.76}  & 625.42  & 98.97  & 616.24  & 623.58  & 100.81 & 615.74  & 623.51  & 99.11   & 0 & 0.2\\
        r7-70  &766.05 & 784.54 & 273.52 & 763.46  & 785.69  & 275.48 & \textbf{760.09}  & 783.13  & 280.49 & 761.58  & 778.04  & 273.50  &1 &1.5 \\
        r7-84  &932.12 & NA    & 584.26 & 897.50  & 932.05  & 488.49 & 897.34  & 915.24  & 446.76 & \textbf{896.91}  & 916.23  & 456.77  & 3 & 3.4\\
        r8-64  &638.36 & 652.30 & 641.63 & 642.34  & 652.65  & 646.45 & 639.01  & 652.80  & 671.52 & \textbf{637.84}  & 652.17  & 719.50  & 0 & 0.2\\
        r8-80  &811.19 & 833.05 & 448.14 & 816.17  & 834.80  & 438.40 & \textbf{808.14}  & 828.89  & 420.03 & 813.16  & 829.92  & 450.94  & 1 & 1.1\\
        r8-96  &NA   & NA    & 617.17 & 1089.18 & 1129.20 & 588.26 & 1060.48 & 1098.13 & 545.21 & \textbf{1058.41} & 1090.04 & 564.49 & 5 &4.6 \\
        \hline
        Avg  & &  &373.89   &  &  &361.23 & &  & 352.29  &  &  &362.17 & 1.4 &1.9 \\
        \hline
        \textbf{$\gamma = 0.7$} & $BC$ & $AC$ &CPU & $BC_2$ & $AC_2$ &CPU$_2$ & $BC_3$ & $AC_3$ &CPU$_3$ & $BC_{\infty}$ & $AC_{\infty}$ &CPU$_{\infty}$ &$N_{max}^s$ &$\overline{N_{avg}^s}$\\
        \hline
        r5-60  &NA & NA & 507.76  & 731.84 & 770.95 & 484.01  & \textbf{704.97}  & 725.74  & 483.86  & 708.54  & 723.73  & 492.51  &5 &6.9 \\
        r6-48  &NA & NA & 502.21  & 518.87 & 540.88 & 507.06  & 509.80  & 525.98  & 486.31  & \textbf{509.76}  & 525.10  & 483.94  & 3 &5.0 \\
        r6-60  &NA & NA & 327.25  & 716.48 & 741.76 & 300.67  & 700.82  & 713.33  & 306.60  & \textbf{697.57}  & 711.52  & 289.76  &6 &7.0 \\
        r6-72  &NA & NA & 590.56  & 920.61 & NA    & 605.16  & 798.26  & 817.20  & 561.24  & \textbf{796.19}  & 826.48  & 574.02   &4 &8.4 \\
        r7-56  &NA & NA & 221.09  & 644.19 & 662.06 & 208.57  & \textbf{622.66}  & 640.69  & 210.29  & 625.91  & 641.82  & 212.05  & 3 &7.0 \\
        r7-70  &NA & NA & 510.60  & 866.06 & NA    & 507.14  & \textbf{777.85}  & 803.20  & 465.43  & 781.56  & 800.35  & 480.03   & 6 &7.8 \\
        r7-84  &NA & NA & 790.95  & NA    & NA    & 753.17  & \textbf{906.14}  & 938.15  & 623.70  & 915.61  & 938.49  & 705.25   & 6 &8.7 \\
        r8-64  &NA & NA & 1207.35 & 664.02 & 698.61 & 1170.20 & \textbf{647.02}  & 666.20  & 1185.16 & 649.93  & 668.48  & 1290.02  &7 &6.3 \\
        r8-80  &NA & NA & 868.04  & 966.47 & NA    & 846.51  & \textbf{829.54}  & 857.56  & 707.30  & 843.26  & 865.90  & 744.33 & 4 &8.2 \\
        r8-96  &NA & NA & 860.97  & NA    & NA    & 845.14  & 1105.82 & 1145.82 & 646.04  & \textbf{1097.76} & 1136.43 & 806.99  &7 &11.2 \\
        \hline
        Avg  & &  &638.68   &  &  &622.76 & &  &567.59   &  &  &607.89  & 5.1 &7.7 \\
        \hline
    \end{tabular}
    \end{center}
\end{table}

\end{APPENDICES}

\end{document}